\documentclass[a4paper]{article}
\usepackage{amsmath}\usepackage{epsf,amsfonts,amsthm}
\usepackage{fontenc,indentfirst, delarray,amsfonts,amsmath,amssymb}

\newcommand{\be}{\begin{equation}}
\newcommand{\ee}{\end{equation}}
\newcommand{\bea}{\begin{eqnarray*}}
\newcommand{\eea}{\end{eqnarray*}}
\newcommand{\w}{\wedge}
\newcommand{\p}{\partial}
\newcommand{\la}{\langle}
\newcommand{\ra}{\rangle}
\newcommand{\raa}{\rightarrow}
\newcommand{\Ci}{C^{\infty}}
\newcommand{\E}{\ell}
\newcommand{\N}{\mathbb{N}}

\newcommand{\R}{\mathbb{R}}

\newcommand{\C}{\mathbb{C}}

\newcommand{\F}{\mathbb{F}}
\newcommand{\n}{\nabla}

\newcommand{\D}{(x_1^2+x_2^2)x_3}

\newcommand{\lp}{\left(}
\newcommand{\rp}{\right)}

\newcommand{\op}[1]{\!\!\mathop{\rm ~#1}\nolimits}

\mathchardef\za="710B  
\mathchardef\zb="710C  
\mathchardef\zg="710D  
\mathchardef\zd="710E  
\mathchardef\zve="710F 
\mathchardef\zz="7110  
\mathchardef\zh="7111  
\mathchardef\zy="7112 
\mathchardef\zi="7113  
\mathchardef\zk="7114  
\mathchardef\zl="7115  
\mathchardef\zm="7116  
\mathchardef\zn="7117  
\mathchardef\zx="7118  
\mathchardef\zp="7119  
\mathchardef\zr="711A  
\mathchardef\zs="711B  
\mathchardef\zt="711C  
\mathchardef\zu="711D  
\mathchardef\zf="711E 
\mathchardef\zq="711F  
\mathchardef\zc="7120  
\mathchardef\zw="7121  
\mathchardef\ze="7122  
\mathchardef\zvy="7123  
\mathchardef\zvw="7124  
\mathchardef\zvr="7125 
\mathchardef\zvs="7126 
\mathchardef\zvf="7127  
\mathchardef\zG="7000  
\mathchardef\zD="7001  
\mathchardef\zY="7002  
\mathchardef\zL="7003  
\mathchardef\zX="7004  
\mathchardef\zP="7005  
\mathchardef\zS="7006  
\mathchardef\zU="7007  
\mathchardef\zF="7008  
\mathchardef\zW="700A  

\newtheorem{theo}{Theorem}
\newtheorem{prop}{Proposition}
\newtheorem{lem}{Lemma}
\newtheorem{cor}{Corollary}

\newtheorem{defi}{Definition}

\begin{document}

\title{Strongly $r$-matrix induced tensors,\\ Koszul cohomology, and\\ arbitrary-dimensional quadratic Poisson cohomology}
\author{Mourad Ammar \footnote{University of Luxembourg, Campus Limpertsberg, Institute
of Mathematics, 162A, avenue de la Fa\"iencerie, L-1511 Luxembourg
City, Grand-Duchy of Luxembourg, E-mail: mourad.ammar@uni.lu,
guy.kass@uni.lu, norbert.poncin@uni.lu. The research of M. Ammar
and N. Poncin was supported by grant R1F105L10. The last author
also thanks the Erwin Schr\"odinger Institute in Vienna for
hospitality and support during his visits in 2006 and 2007.} , Guy
Kass $^*$, Mohsen Masmoudi \footnote{Universit\'e Henri
Poincar\'e, Institut Elie Cartan, B.P. 239, F-54 506
Vandoeuvre-les-Nancy Cedex, France, E-mail:
Mohsen.Masmoudi@iecn.u-nancy.fr} , Norbert Poncin $^*$}\maketitle

\begin{abstract}
We introduce the concept of strongly $r$-matrix induced ({\small
SRMI}) Poisson structure, report on the relation of this property
with the stabilizer dimension of the considered quadratic Poisson
tensor, and classify the Poisson structures of the Dufour-Haraki
classification (DHC) according to their membership of the family
of {\small SRMI} tensors. One of the main results of our work is a
generic cohomological procedure for {\small SRMI} Poisson
structures in arbitrary dimension. This approach allows
decomposing Poisson cohomology into, basically, a Koszul
cohomology and a relative cohomology. Moreover, we investigate
this associated Koszul cohomology, highlight its tight connections
with Spectral Theory, and reduce the computation of this main
building block of Poisson cohomology to a problem of linear
algebra. We apply these upshots to two structures of the DHC and
provide an exhaustive description of their cohomology. We thus
complete our list of data obtained in previous works, see
\cite{MP} and \cite{AMPN}, and gain fairly
good insight into the structure of Poisson cohomology.\\



\end{abstract}

\section{Introduction}

In a graded Lie algebra (gLa) $({\cal L},[.,.])$, ${\cal
L}=\oplus_i{\cal L}^{i},$ any element with degree $1$ that squares
to $0$, generates a differential graded Lie algebra (dgLa) $({\cal
L},[.,.],\p_{\zL})$, $\p_{\zL}:=[\zL,.]$, and a gLa $H({\cal
L},[.,.],\p_{\zL})$ in cohomology. It is interesting to note that,
depending on the initial algebra, such a 2-nilpotent degree $1$
element is e.g. an associative algebra structure, a Lie algebra
structure, or a Poisson structure, and the associated cohomology
is the adjoint Hochschild, the adjoint Chevalley-Eilenberg, and
the Lichnerowicz-Poisson (LP) (or simply Poisson) cohomology,
respectively. Let us recall that the LP-dgLa is implemented by the
shifted Grassmann algebra $({\cal X}(M)[1],\w,[.,.]_{\op{SN}})$,
${\cal X}(M)=\zG(\w TM)$, of polyvectors of a manifold $M$,
endowed with the Schouten-Nijenhuis bracket $[.,.]_{\op{SN}}$
(whereas the Hochschild (resp. the Chevalley-Eilenberg) dgLa is
generated by the space of multilinear (resp. skew-symmetric
multilinear) mappings of the underlying vector space, endowed with
the Gerstenhaber (resp. the Nijenhuis-Richardson) graded Lie
bracket).

Alternatively, LP-cohomology can be viewed as the Lie algebroid
(Lad) cohomology of the Lie algebroid $(T^*M,\{.,.\},\sharp)$
canonically associated with an arbitrary Poisson manifold
$(M,\zL)$ (usual notations). The cohomology of a Lad $(E\raa
M,[\![.,.]\!],\zr)$ (or, equivalently, a $Q$-structure on a
supermanifold) is defined as the cohomology of the
Chevalley-Eilenberg subcomplex of the representation
$\zr:\zG(E)\raa \op{Der}(\Ci(M))$, made up by tensorial cochains.
Algebraically, LP-cohomology is defined as the adjoint
Chevalley-Eilenberg cohomology of any Poisson-Lie algebra,
restricted to the cochain subspace of skew-symmetric
multiderivations.

More detailed descriptions of Poisson cohomology can be found e.g. in \cite{AL} or \cite{IV}.\\

Many papers on Poisson cohomology and Poisson homology,
\cite{ko:homology}, \cite{br:homology}, have been published during
the last decades. Cohomology of regular Poisson manifolds,
\cite{IV3}, \cite{PX}, (co)homology and resolutions,
\cite{hu:resolutions}, duality, \cite{hu:duality},
\cite{xu:duality}, \cite{ELW}, cohomology in low dimensions or
specific cases, \cite{NN}, \cite{VLG}, \cite{AnGa}, \cite{PM1},
\cite{PM3}, \cite{RV}, \cite{DR}, \cite{AP}, extensions of Poisson
cohomology, e.g. Lie algebroid cohomology, Jacobi cohomology,
Nambu-Poisson cohomology, double Poisson cohomology, and graded
Jacobi cohomology, \cite{LMP}, \cite{ILLMP}, \cite{PM2},
\cite{GM}, \cite{LLMP}, \cite{NN2}, \cite{PW}, are only some of
the investigated problems. Let us also mention our own works,
\cite{MP}, \cite{AMPN}, in which we suggest an approach to the
cohomology of the Poisson tensors of the
Dufour-Haraki classification ({\small DHC}).\\

In this paper, we focus on the formal LP-cohomology associated
with the quadratic Poisson tensors ({\small QPT}) $\zL$ of $\R^n$
that read as real linear combination
\begin{equation}\label{admissible} \zL=\sum_{i<j}\alpha^{ij}Y_i\w
Y_j=:\sum_{i<j}\alpha^{ij}Y_{ij},\quad \alpha^{ij}\in\R
\end{equation} of the wedge products of $n$ commuting linear
vector fields $Y_1,\ldots,Y_n$, such that $Y_1\w\ldots\w
Y_n=:Y_{1\ldots n}\neq0$. Let us recall that ``formal'' means that
we substitute the space $\R[[x_1,\ldots,x_n]]\otimes \w\R^n$ of
multivectors with coefficients in the formal series for the usual
Poisson cochain space ${\cal X}(\R^n)=\Ci(\R^n)\otimes\w\R^n$.
Furthermore, the reader may think about {\small QPT} of type
(\ref{admissible}) as {\small QPT} implemented by a classical
$r$-matrix in their stabilizer for the canonical matrix action.

Hence, in Section 2, we are interested in the characterization of
the {\small QPT} that are images of a classical $r$-matrix. We
comment on the tight relation between the fact that a {\small QPT}
is induced by an $r$-matrix and the dimension of its stabilizer.
We prove that if the stabilizer of a given {\small QPT} $\zL$ of
$\R^n$ contains $n$ commuting linear vector fields $Y_i$, such
that $Y_{1...n}\neq 0$, then $\zL$ is implemented by an $r$-matrix
in its stabilizer, see Theorem \ref{1}. In the following, we refer
to such tensors as strongly $r$-matrix induced ({\small SRMI})
structures and show that any structure of the {\small DHC}
decomposes into the sum of a major {\small SRMI} structure and a
small compatible (mostly exact) Poisson tensor, see Theorem
\ref{ClassTheo}. This decomposition constitutes the foundation of
our cohomological techniques proposed in \cite{MP} and
\cite{AMPN}. The preceding description and the philosophy of the
mentioned cohomological modus operandi allow understanding that
our splitting is in some sense in opposition to the one proven in
\cite{LX} that incorporates the largest possible part of the
Poisson tensor into the exact term.

In \cite{MP}, two of us developed a cohomological method in the
Euclidean Three-Space that led to a significant simplification of
LP-cohomology computations for the {\small SRMI} structures of the
{\small DHC}. Section 3 of the present note aims for extension of
this procedure to arbitrary dimensional vector spaces. Nontrivial
lemmata allow injecting the space ${\cal R}$ of ``real''
LP-cochains (formal multivector fields) into a larger space ${\cal
P}$ of ``potential'' cochains, see Theorem \ref{injection}, and
identifying the natural extension to ${\cal P}$ of the
LP-differential as the Koszul differential associated with $n$
commuting endomorphisms $X_i-(\op{div}X_i)\op{id}$,
$X_i=\sum_j\alpha^{ij}Y_j$, $\za^{ji}=-\za^{ij}$, of the space
made up by the polynomials on $\R^n$ with some fixed homogeneous
degree, Theorems \ref{PoissKoszCob} and \ref{PoissKosz}. We then
choose a space ${\cal S}$ supplementary to ${\cal R}$ in ${\cal
P}$ and show that the LP-differential induces a differential on
${\cal S}$. Eventually, we end up with a short exact sequence of
differential spaces and an exact triangle in cohomology. It could
be proven that the LP-cohomology (${\cal R}$-cohomology) reduces,
essentially, to the above-depicted Koszul cohomology (${\cal
P}$-cohomology) and a relative cohomology (${\cal S}$-cohomology),
see Theorem \ref{PoissKoszRel}.

In order to take advantage of these upshots, we investigate in
Section 4 the Koszul cohomology associated to $n$ commuting linear
operators on a finite-dimensional complex vector space. We prove a
homotopy-type formula, see Proposition \ref{HomotKosz},
and---using spectral properties---we show that the Koszul
cohomology is, roughly spoken, located inside (a direct sum of
intersections of) the kernels of some transformations that can be
constructed recursively from the initially considered operators,
Proposition \ref{KoszCohoBasic} and Corollary \ref{KoszCohoFin}.

In Section 5, we apply this result, gain valuable insight into the
structure of the Koszul cohomology implemented by {\small SRMI}
tensors, and show that in order to compute this central part of
Poisson cohomology it basically suffices to solve triangular
systems of linear equations.

Section 6 contains a full description of the LP-cohomology spaces
of structures $\zL_3$ and $\zL_9$ of the {\small DHC}.

Eventually, the aforementioned general upshots and our growing
list of explicit data allow describing the main LP-cohomological
phenomena, see Section 7.

\section{Characterization of strongly $r$-matrix induced Poisson structures}\label{admstr}

\subsection{Stabilizer dimension and $r$-matrix generation}

Poisson structures implemented by an $r$-matrix are of importance,
e.g. in Deformation Quantization, especially in view of Drinfeld's
method. In the following, we report on an idea regarding
generation of quadratic Poisson tensors by classical
$r$-matrices.\\

Set $G=\op{GL}(n,\R)$ and $\frak{g}=\op{gl}(n,\R)$. The Lie
algebra isomorphism between $\frak{g}$ and the algebra ${\cal
X}^1_0(\R^n)$ of linear vector fields, extends to a Grassmann
algebra and a graded Poisson-Lie algebra homomorphism
$J:\wedge\frak{g}\raa \oplus_k\lp{\cal
S}^k\R^{n*}\otimes\w^k\R^n\rp.$ It is known that its restriction
$$J^k:\w^k\frak{g}\raa {\cal S}^k\R^{n*}\otimes\w^k\R^n$$ is onto,
but has a non-trivial kernel if $k,n\ge 2$. In particular,
\[J^3[r,r]_{\op{SN}}=[J^2r,J^2r]_{\op{SN}},\;r\in \frak{g}\wedge \frak{g},\] where
$[.,.]_{\op{SN}}$ is the Schouten-Nijenhuis bracket. These
observations allow to understand that the characterization of the
quadratic Poisson structures that are implemented by a classical
$r$-matrix, i.e. a bimatrix $r\in\frak{g}\w\frak{g}$ that verifies
the Classical Yang-Baxter Equation $[r,r]_{\op{SN}}=0$, is an open
problem.

Quadratic Poisson tensors $\zL_1$ and $\zL_2$ are equivalent if
and only if there is $A\in G$ such that $A_*\zL_1=\zL_2,$ where
$*$ denotes the standard action of $G$ on tensors of $\R^n$. As
$J^2$ is a $G$-module homomorphism, i.e.
$$A_*(J^2r)=J^2(\op{Ad}(A)r),\;A\in G, r\in\frak{g}\w\frak{g},$$ the $G$-orbit
of a quadratic Poisson structure $\zL=J^2r$ is the pointwise
$J^2$-image of the $G$-orbit of $r$. Furthermore, representation
$\op{Ad}$ acts by graded Lie algebra homomorphisms, i.e.
$$\op{Ad}(A)[r,r]_{\op{SN}}=[\op{Ad}(A)r,\op{Ad}(A)r]_{\op{SN}}.$$
Hence, if $\zL=J^2r$, where $r$ is a classical $r$-matrix, the
whole orbit of this quadratic Poisson tensor is made up by
$r$-matrix induced structures.

Of course, any quadratic Poisson tensor $\zL$ is implemented by
bimatrices $r\in\frak{g}\w\frak{g}$. In order to determine wether
the $G$-orbit $O_{\zL}$ of this tensor is generated by
$r$-matrices, we have to take an interest in the preimage
$$(J^2)^{-1}(O_{\zL})=\cup_{r\in(J^2)^{-1}\zL}O_r,$$
composed of the $G$-orbits $O_r$ of all the bimatrices $r$ that
are mapped on $\zL$ by $J^2$. We claim that the chances that a
fiber of this bundle is located inside $r$-matrices are the
bigger, the smaller is $O_{\zL}$. In other words, the dimension of
the isotropy Lie group $G_{\zL}$ of $\zL$, or of its Lie algebra,
the stabilizer
$$\frak{g}_{\zL}=\{a\in \frak{g}: [\zL,Ja]_{\op{SN}}=0\}$$ of $\zL$ for the
corresponding infinitesimal action, should be big enough. In
addition to the ostensible intuitive clearness of this conjecture,
positive evidence comes from the fact that, in $\R^3$, the Poisson
tensor $\zL=(x_1^2+x_2x_3)\p_{23}$, $\p_{23}:=\p_2\w\p_3$,
$\p_i:=\p/\p_{x_i}$, is not $r$-matrix induced, see \cite{MMR},
and the dimension of its stabilizer is $\op{dim}\frak{g}_{\zL}=2$,
as well as from the following theorem (we implicitly identify
stabilizer $\frak{g}_{\zL}\subset {\frak g}$ and the (isomorphic)
Lie subalgebra $J^1\frak{g}_{\zL}=\{Y\in{\cal X}^1_0(\R^n):
[\zL,Y]_{\op{SN}}=0\}\subset {\cal X}^1_0(\R^n)$ of linear vector
fields of $\R^n$).

\begin{theo}\label{1}Let $\zL$ be a quadratic Poisson tensor of $\R^n$. If its stabilizer
${\frak g}_{\zL}$ contains $n$ commuting linear vector fields
$Y_i$, $i\in\{1,\ldots,n\}$, such that $Y_1\w\ldots\w Y_n\neq0$,
then $\Lambda$ is implemented by a classical $r$-matrix that
belongs to the stabilizer, i.e. $\zL=J^2a,$ $[a,a]_{\op{SN}}=0$,
$a\in\frak{g}_{\zL}\w\frak{g}_{\zL}$.
\end{theo}
{\it Proof}. Let $(x_1,\ldots,x_n)$ be the canonical coordinates
of $\R^n$. Set $\p_r=\p_{x_r}$ and $Y_i=\sum_{r=1}^n\E_{ir}\p_r$,
with $\E\in\op{gl(n,\R^{n*})}$. The determinant $D=\op{det}\E$
does not vanish everywhere, since $Y_{1\ldots n}=D\,\p_{1\ldots
n}$ and $Y_{1\ldots n}\neq 0$. At any point of the nonempty open
subset $Z=\{x\in\R^n, D(x)\neq 0\}$ of $\R^n$, the $Y_i$ form a
basis of the corresponding tangent space of $\R^n.$ Moreover, in
$Z$, we get
$$\p_{ij}=D^{-2}\sum_{k<l}\left(\mathbf{L}^k_{\;i}\mathbf{L}^l_{\;j}-\mathbf{L}^l_{\;i}\mathbf{L}^k_{\;j}\right)Y_{kl}
=:D^{-2}\sum_{k<l}Q_{ij}^{kl}Y_{kl},$$ where $\mathbf{L}$ denotes
the matrix of maximal algebraic minors of $\E$, and where
$Q^{kl}_{ij}\in{\cal S}^{2n-2}\R^{n*}$. Hence, if the quadratic
Poisson tensor $\zL$ reads $\zL=\sum_{i<j}\zL^{ij}\p_{ij}$,
$\zL^{ij}\in{\cal S}^2\R^{n*}$, we have in $Z$,
$$\zL=D^{-2}\sum_{k<l}\sum_{i<j}\zL^{ij}Q_{ij}^{kl}Y_{kl}=:D^{-2}\sum_{k<l}P^{kl}Y_{kl},$$ where $P^{kl}\in{\cal S}^{2n}\R^{n*}$.
We now prove that the rational functions $D^{-2}P^{kl}$ are
actually constants. Since the $Y_i$ are commuting vector fields in
${\frak g}_{\zL}$, the commutation relations $[Y_i,Y_j]=0$ and
$[\zL,Y_i]_{\op{SN}}=0$, $i,j\in\{1,\ldots,n\}$, hold true. It
follows that
$$Y_i\left(D^{-2}P^{kl}\right)=\sum_{r=1}^n\E_{ir}\p_r\left(D^{-2}P^{kl}\right)=0,\quad i\in\{1,\ldots,n\},$$
everywhere in $Z$, and, as $\E$ is invertible in $Z$, that
$\p_r\left(D^{-2}P^{kl}\right)=0,\,r\in\{1,\ldots,n\}.$ Hence,
$P^{kl}=\za^{kl}D^2,\,\za^{kl}\in\R,$ in each connected component
of $Z$. As these components are open subsets of $\R^n$, the last
result holds in $\R^n$ (in particular the constants $\za^{kl}$
associated with different connected components coincide).
Eventually,
$$\zL=\sum_{k<l}\za^{kl}Y_{kl}=J^2\left(\sum_{k<l}\za^{kl}a_{kl}\right),$$
where $a_i=(J^1)^{-1}Y_i\in{\frak g}_{\zL}$. Since the $a_i$ are
(just as the $Y_i$) mutually commuting, it is clear that the
bimatrix $r=\sum_{k<l}\za^{kl}a_{kl}\in{\frak g}_{\zL}\w{\frak
g}_{\zL}$ verifies the classical Yang-Baxter
equation.\quad\rule{1.5mm}{2.5mm}

\begin{defi} We refer to a quadratic Poisson structure $\zL$ that is implemented by a classical $r$-matrix
$r\in\frak{g}_{\zL}\w\frak{g}_{\zL}$, where ${\frak g}_{\zL}$
denotes the stabilizer of $\zL$ for the canonical matrix action,
as a strongly $r$-matrix induced $(${\small SRMI }$)$ tensor.
\end{defi}

\subsection{Classification theorem in Euclidean Three-Space}

Two concepts of exact Poisson structure---tightly related with two
special cohomology classes---are used below. Let $\zL$ be a
Poisson tensor on a smooth manifold $M$ oriented by a volume
element $\zW$. We say that $\zL$, which is of course a
LP-$2$-cocycle, is LP-exact (Lichnerowicz-Poisson), if
$$\zL=[\zL,X]_{\op{SN}},\;X\in{\cal X}^1(M),$$ [vector field $X$ is
called Liouville vector field and the cohomology class of $\zL$ is
the obstruction to infinitesimal rescaling of $\zL$], and we term
$\zL$ K-exact (Koszul), if
$$\zL=\delta(T),\;T\in{\cal X}^3(M).$$ Operator $\delta:=\zf^{-1}\circ
d\circ\zf$ is the pullback of the de Rham differential $d$ by the
canonical vector space isomorphism $\zf:=i_{\cdot}\zW$. Although
introduced earlier, the generalized divergence $\delta$
($\delta(X)=\op{div}_{\zW}X,\;X\in{\cal X}^1(M)$) is prevalently
attributed to J.-L. Koszul. The curl vector field
$K(\zL):=\delta(\zL)$ of $\zL$ (if $\zW$ is the standard volume of
$\R^3$ and $\zL$ is identified with a vector field $\vec{\zL}$ of
$\R^3$, $K(\zL)$ coincides with the standard curl
$\vec{\n}\w\vec{\zL}$) is an LP-$1$-cocycle (which maps a function
to the divergence of its Hamiltonian vector field, and the
cohomology class of which is the well-known modular class of $\zL$
[this class is independent of $\zW$ and is the obstruction to
existence on $M$ of a measure preserved by all Poisson
automorphisms] that is relevant e.g. in the classification of
Poisson structures, see \cite{DH}, \cite{GMP}, \cite{LX}, and in
Poincar\'e duality, see \cite{ELW}, \cite{ILLMP}). In $\R^n$,
$n\ge 3$, a Poisson tensor $\zL$ is K-exact, if and only if it is
``irrotational'', i.e. $K(\zL)=0$, and in $\R^3$, K-exact means
``function-induced'', i.e.
$$\zL={
\zP_f}:=\p_1f\,\p_{23}+\p_2f\,\p_{31}+\p_3f\,\p_{12},\;f\in\Ci(\R^3).$$
The K-exact quadratic Poisson tensors $\zP_p$ of $\R^3$, i.e. the
K-exact Poisson structures that are induced by a homogeneous
polynomial $p\in{\cal S}^3\R^{3*}$, represent class 14 of the DHC.
The cohomology of this class has been studied in \cite{AP}
(actually the author deals with structures $\zP_p$ implemented by
a weight homogeneous polynomial $p$ with an isolated singularity).
Hence, class 14 of the DHC will not be examined in the current
work.

Let us also recall that two Poisson tensors $\zL_1$ and $\zL_2$
are compatible, if their sum is again a Poisson structure, i.e. if
$[\zL_1,\zL_2]_{\op{SN}}=0$.

The following theorem classifies the quadratic Poisson classes
according to their membership of the family of strongly $r$-matrix
induced structures. Furthermore, we show that any structure reads
as the sum of a {\it major} strongly $r$-matrix induced tensor and
a {\it small} compatible Poisson structure. On one hand, this
membership entails accessibility to the cohomological technique
exemplified in \cite{MP}, on the other, this splitting---which, by
the way, differs from the decomposition suggested in \cite{LX} in
the sense that we incorporate the biggest possible part of the
structure into the strongly induced term---is of particular
importance with regard to the cohomological approach detailed in
\cite{AMPN}.
\begin{theo}\label{ClassTheo}
Let $a,b,c\in\R$ and let $\zL_i$ $(i\in\{1,\ldots,13\})$ be the
quadratic Poisson tensors of the DHC, see \cite{DH}. Denote the
canonical coordinates of $\R^3$ by $x,y,z$ (or $x_1,x_2,x_3$) and
the partial derivatives with respect to these coordinates by
$\p_1,\p_2,\p_3$ $(\p_{ij}=\p_i\w\p_j)$.

If $\op{dim}\frak{g}_{\zL}>3$ (subscript $i$ omitted), there are
mutually commuting linear vector fields $Y_1,Y_2,Y_3$, such that
$$\zL=\za Y_{23}+\zb Y_{31}+\zg Y_{12}\;\;(\za,\zb,\zg\in\R),$$ so that
$\zL$ is strongly $r$-matrix induced ({\small SRMI}), i.e.
implemented by a classical $r$-matrix in
$\frak{g}_{\zL}\w\frak{g}_{\zL}$. In the following classification
of the quadratic Poisson tensors with regard to property {\small
SRMI}, we decompose each not {\small SRMI} tensor into the sum of
a major {\small SRMI} structure and a smaller compatible quadratic
Poisson tensor.

\begin{itemize}

\item Set $Y_1=x\p_1,Y_2=y\p_2,Y_3=z\p_3$

\begin{enumerate}

\item $\zL_1=a\,yz\p_{23}+b\,xz\p_{31}+c\,xy\p_{12}$ is {\small
SRMI} for all values of the parameters $a,b,c$. More precisely,
\[\zL_1=a\,Y_{23}+b\,Y_{31}+c\,Y_{12}\]

\item $\zL_4=a\,yz\p_{23}+a\,xz\p_{31}+(b\,xy+z^2)\p_{12}$ is not
{\small SRMI} if and only if $(a,b)\neq (0,0)$. We have,
$$\zL_4=a(Y_{23}+Y_{31})+b\,Y_{12}+{\frac{1}{3}}\zP_{z^3}$$

\end{enumerate}

\item Set $Y_1=x\p_1+y\p_2,Y_2=x\p_2-y\p_1,Y_3=z\p_3$

\begin{enumerate}

\item $\zL_2=(2ax-by)z\p_{23}+(bx+2ay)z\p_{31}+a(x^2+y^2)\p_{12}$
is {\small SRMI} for any $a,b$. More precisely,
$$\Lambda_2=2a\,Y_{23}+b\,Y_{31}+a\,Y_{12}$$

\item $\zL_7=\lp(2a+c)x-by\rp z\p_{23}+\lp bx+(2a+c)y\rp
z\p_{31}+a(x^2+y^2)\p_{12}$ is {\small SRMI} for all $a,b,c$. More
precisely,
$$\zL_7=(2a+c)Y_{23}+b\,Y_{31}+a\,Y_{12}$$

\item
$\zL_8=a\,xz\p_{23}+a\,yz\p_{31}+\lp\frac{a+b}{2}(x^2+y^2)\pm
z^2\rp\p_{12}$ is not {\small SRMI} if and only if $(a,b)\neq
(0,0)$. We have,
$$\zL_8=a\,Y_{23}+\frac{a+b}{2}Y_{12}\pm\frac{1}{3}\zP_{z^3}$$

\end{enumerate}

\item Set $Y_1=x\p_1+y\p_2,Y_2=x\p_2,Y_3=z\p_3$

\begin{enumerate}

\item $\zL_3=(2x-a\,y)z\p_{23}+a\,xz\p_{31}+x^2\p_{12}$ is {\small
SRMI} for any $a$. More precisely,
$$\zL_3=2Y_{23}+a\,Y_{31}+Y_{12}$$

\item $\zL_5=\lp(2a+1)x+y\rp z\p_{23}-xz\p_{31}+a\,x^2\p_{12}$
$(a\neq -\frac{1}{2})$ is {\small SRMI} for any $a$. More
precisely,
$$\zL_5=(2a+1)Y_{23}-Y_{31}+a\,Y_{12}$$

\item $\zL_6=a\,yz\p_{23}-a\,xz\p_{31}-\frac{1}{2}x^2\p_{12}$ is
{SRMI} for any $a$. More precisely,
$$\zL_6=-a\,Y_{31}-\frac{1}{2}Y_{12}$$

\end{enumerate}

\item Set $Y_1={\cal
E}:=x\p_1+y\p_2+z\p_3,Y_2=x\p_2+y\p_3,Y_3=x\p_3$

\begin{enumerate}

\item
$\zL_9=(ax^2-\frac{1}{3}y^2+\frac{1}{3}xz)\p_{23}+\frac{1}{3}xy\p_{31}-\frac{1}{3}x^2\p_{12}$
is {\small SRMI} for any $a$. More precisely,
$$\zL_9=a\,Y_{23}-\frac{1}{3}Y_{12}$$

\item $\zL_{10}=\lp
a\,y^2-(4a+1)xz\rp\p_{23}+(2a+1)xy\p_{31}-(2a+1)x^2\p_{12}$ is not
{\small SRMI} if and only if $a\neq-\frac{1}{3}$. We have,
$$\zL_{10}=-(2a+1)Y_{12}+(3a+1)(y^2-2xz)\p_{23}$$

\end{enumerate}

\item Set $Y_1={\cal E},Y_2=x\p_2,Y_3=\lp
a\,x+(3b+1)z\rp\p_3$

\begin{enumerate}

\item $\zL_{11}=\lp
a\,x^2+(2b+1)xz\rp\p_{23}+(b\,x^2+c\,z^2)\p_{12}$ $(a=0)$ is not
{\small SRMI} if and only if $c\neq 0$. We have,
$$\zL_{11}=Y_{23}+b\,Y_{12}+\frac{c}{3}\zP_{z^3}$$

\item $\zL_{12}=\lp
a\,x^2+(2b+1)xz\rp\p_{23}+(b\,x^2+c\,z^2)\p_{12}$ $(a=1)$ is not
{\small SRMI} if and only if $c\neq 0$. We have,
$$\zL_{12}=Y_{23}+b\,Y_{12}+\frac{c}{3}\zP_{z^3}$$

\item $\zL_{13}=\lp
a\,x^2+(2b+1)xz+z^2\rp\p_{23}+(b\,x^2+c\,z^2+2xz)\p_{12}$ is not
{\small SRMI} for any $a,b,c$. We have,
$$\zL_{13}=Y_{23}+b\,Y_{12}+\zP_{\frac{c}{3}z^3+xz^2}$$

\end{enumerate}
\end{itemize}
\end{theo}

{\it Proof.} Let us first mention that the specified basic fields
$Y_1,Y_2,Y_3$ have been read in the stabilizers of the considered
Poisson tensors, but that we refrain from publishing the often
fairly protracted stabilizer-computations. Indeed, once the vector
fields $Y_i$ are known, it is easily checked that, in the {\small
SRMI} cases, they verify the assumptions of Theorem \ref{1}. Thus
the corresponding Poisson structures are actually {\small SRMI}
tensors. In order to ascertain that a quadratic Poisson structure
$\zL$ is not {\small SRMI}, it suffices to prove that $\zL\notin
J^2(\frak{g}_{\zL}\w\frak{g}_{\zL})$. This will be done
thereinafter. All the quoted decompositions can be directly
verified. In most instances, the twist is obviously Poisson, so
that compatibility follows. In the case of $\zL_{10}$, the twist
$\zL_{10,\op{II}}=(y^2-2xz)\p_{23}$ is a non-K-exact Poisson
structure. This is a direct consequence of the result
$K(\zL_{10,\op{II}})=\vec{\n}\w\vec{\zL}_{10,\op{II}}=-2x\p_2-2y\p_3\neq
0$ and the handy formula
$$[P,Q]_{\op{SN}}=(-1)^pD(P\w Q)-D(P)\w Q-(-1)^pP\w D(Q),\forall P\in{\cal X}^p(M),Q\in{\cal
X}^q(M).$$ The statement regarding the dimension of stabilizer
$\frak{g}_{\zL}$ is obvious in view of the following main part of
this proof.

Denote by $E_{ij}$ ($i,j\in\{1,2,3\}$) the canonical basis of
$\op{gl}(3,\R)$.
\begin{itemize}

\item For $\zL_4$, if $(a,b)\neq (0,0)$, stabilizer
$\frak{g}_{\zL_4}$ and the image $J^2(\frak{g}_{\zL_4}\w
\frak{g}_{\zL_4})$ are generated by
\[(\frac{1}{2}E_{11}+E_{22},\frac{1}{2}E_{11}+E_{33})\;\;\;\mbox{and}\;\;\;yz\p_{23}-\frac{1}{2}xz\p_{31}-\frac{1}{2}xy\p_{12},\]
respectively. Hence, $\zL_4$ is not {\small SRMI}.

\item For $\zL_8$, if $(a,b)\neq (0,0)$, the generators of
$\frak{g}_{\zL_8}$ and $J^2(\frak{g}_{\zL_8}\w \frak{g}_{\zL_8})$
are
\[(E_{11}+E_{22}+E_{33},E_{12}-E_{21})\;\;\;\mbox{and}\;\;\;-xz\p_{23}-yz\p_{31}+(x^2+y^2)\p_{12}.\] So
$\zL_8$ is not {\small SRMI}.

\item For $\zL_{10},$ if $a\neq-\frac{1}{3}$, the generators are
\[(E_{11}+E_{22}+E_{33},E_{12}+E_{23})\;\;\;\mbox{and}\;\;\;(y^2-xz)\p_{23}-xy\p_{31}+x^2\p_{12}.\]

\item For $\zL_{11},\,c\neq 0$, $\zL_{12},\,c\neq 0$, and
$\zL_{13}$, the generators are
\[(E_{11}+E_{22}+E_{33},E_{12},E_{32})\;\;\;\mbox{and}\;\;\;(-xz\p_{23}+x^2\p_{12},z^2\p_{23}-xz\p_{12}).\;\rule{1.5mm}{2.5mm}\]
\end{itemize}

{\bf Remarks.} \begin{itemize} \item For $\Lambda=\Lambda_i$,
$i\in\{11,12,13\}$, $c\neq0$ if $i\in\{11,12\}$, the dimension of
the stabilizer is $\op{dim}\frak{g}_{\zL}=3$, whereas
$J\frak{g}_{\zL}\w J\frak{g}_{\zL}\w J\frak{g}_{\zL}=\{0\}.$
Hence, if the dimension of the stabilizer coincides with the
dimension of the space, the Poisson structure is not necessarily a
{\small SRMI} tensor. \item For $\zL_{10}$ e.g., the decomposition
proved in \cite{LX} yields
$$\zL_{10}=-\frac{1}{3}Y_{12}\;+\;\zP_{\frac{c}{3}z^3+xz^2+(b+\frac{1}{3})x^2z+\frac{a}{3}x^3}.$$
\end{itemize}

\section{Poisson cohomology of quadratic structures in a finite-dimensional vector space}

\subsection{Koszul homology and cohomology}

Let $\w=\w_n\la\vec{\zh}\ra$ be the Grassmann algebra on
$n\in\N_0$ generators $\vec{\zh}=(\zh_1,\ldots,\zh_n)$, i.e. the
algebra generated over a field $\mathbb{F}$ of characteristic $0$
(in this work $\F=\R$ or $\F=\C$) by generators
$\zh_1,\ldots,\zh_n$ subject to the anticommutation relations
$\zh_k\zh_{\E}+\zh_{\E}\zh_k=0$, $k,\E\in\{1,\ldots,n\}$. Set
$\w=\oplus_{p=0}^n\w^p$, with obvious notations, and let
$\vec{h}=(h_1,\ldots,h_n)$ be dual generators:
$i_{h_k}\zh_{\E}=\zd_{k\E}$. We also need the creation operator
$e_{\zh_k}:\w\ni\zw\raa\zh_k\,\zw\in\w$ and the annihilation
operator $i_{h_k}:\w\ni\zw\raa i_{h_k}\zw\in\w$, where the
interior product is defined as usual. Eventually, we denote by $E$
a vector space over $\F$ and by $\vec{X}=(X_1,\ldots,X_n)$ an
$n$-tuple of commuting linear operators on $E$.

\begin{defi} The complex $$0\raa E\otimes_{\F} \w^n\raa E\otimes_{\F} \w^{n-1}
\raa\ldots\raa E\otimes_{\F} \w^1\raa E\raa 0,$$ with differential
${\zk}_{\vec{X}}=\sum_{k=1}^n X_k\otimes i_{h_k}$, is the Koszul
chain complex ($K_*$-complex) $K_*(\vec{X},E)$ associated with
$\vec{X}$ on $E$. The Koszul homology group is denoted by
$K\!H_*(\vec{X},E)$.
\end{defi}

\begin{defi} The complex $$0\raa E\raa E\otimes_{\F} \w^{1}
\raa\ldots\raa E\otimes_{\F} \w^{n-1}\raa E\otimes_{\F}\w^{n}\raa
0,$$ with differential ${\cal K}_{\vec{X}}=\sum_{k=1}^n X_k\otimes
e_{\zh_k}$, is the Koszul cochain complex ($K^*$-complex)
$K^*(\vec{X},E)$ associated with $\vec{X}$ on $E$. We denote by
$K\!H^*(\vec{X},E)$ the corresponding Koszul cohomology group.
\end{defi}

Observe that commutation of the $X_k$ and anticommutation of the
$i_{h_k}$ (resp. the $e_{\zh_k}$) entail that $\zk_{\vec{X}}$
(resp. ${\cal K}_{\vec{X}}$) actually squares to $0$.\\

\noindent {\bf Example 1.} It is easily checked that, if we choose
$\F=\R$, $E=\Ci(\R^3)$, $\zh_k=dx_k$ (resp. $\zh_k=\p_k=\p_{x_k}$
and $h_k=dx_k$), and $X_k=\p_k$ ($k\in\{1,2,3\}$, $x_1,x_2,x_3$
canonical coordinates of $\R^3$), the $K^*$-complex (resp. the
$K_*$-complex) is nothing but the de Rham complex $(\zW(\R^3),d)$
(resp. its dual version $({\cal X}(\R^3),\delta)$, see above). Note
that, if we identify the subspaces $\zW^k(\R^3)$ of homogeneous
forms with the corresponding spaces of components $E,E^3,E^3,E$,
this $K^*$-complex reads \begin{equation}0\raa E \stackrel{{\cal
K}=\vec{\n}(.)}{\raa} E^3 \stackrel{{\cal K}=\vec{\n}\w(.)}{\raa}
E^3 \stackrel{{\cal K}=\vec{\n}\cdot(.)}{\raa}E\raa
0,\label{RhamKoszul}\end{equation} with
self-explaining notations.\\

\noindent {\bf Example 2.} For $\F=\R,$ $E={\cal
S}\R^{3*}=\R[x_1,x_2,x_3]$, $\zh_k=\p_k$, $X_k=\frak{m}_{P_k}$
($k\in\{1,2,3\}$, $P_k\in E^{d_k}$, $d_k\in\N$,
$\frak{m}_{P_k}:E\ni Q\raa P_kQ\in E$), the chain spaces of the
$K_*$-complex are the spaces of homogeneous polyvector fields on
$\R^3$ with polynomial coefficients, and anew identification with
the corresponding spaces $E,E^3,E^3,E$ of components, allows to
write this $K_*$-complex in the form \begin{equation}0\raa E
\stackrel{\zk=(.)\vec{P}}{\raa} E^3
\stackrel{\zk=(.)\w\vec{P}}{\raa} E^3
\stackrel{\zk=(.)\cdot\vec{P}}{\raa} E\raa 0,\label{MultiplicKosul}\end{equation} where $\vec{P}=(P_1,P_2,P_3).$\\

\noindent {\bf Remarks.} \begin{itemize}\item Of course, the
Koszul cohomology and homology complexes defined in Example 1 are
exact, expect that $K\!H^0(\vec{\p},\Ci(\R^3))\simeq
K\!H_3(\vec{\p},\Ci(\R^3))\simeq\R$.

\item Let us recall that an $R$--regular sequence on a module $M$
over a commutative unit ring $R$, is a sequence
$(r_1,\ldots,r_d)\in R^d$, such that $r_k$ is not a zero divisor
on the quotient $M/\la r_1, \ldots, r_{k-1}\ra M$,
$k\in\{1,\ldots,d\}$, and $M/\la r_1, \ldots, r_{d}\ra M\neq 0$.
In particular, $x_1,\ldots, x_d$ is a (maximal length) regular
sequence on the polynomial ring $R=\F[x_1,\ldots,x_d]$ (so that
this ring has depth $d$).

It is well-known that the $K_*$-complex described in Example 2 is
exact, except for surjectivity of $\zk=(.)\cdot\vec{P}$, if
sequence $\vec{P}=(P_1,P_2,P_3)$ is regular for $\R[x_1,x_2,x_3]$.
For instance, if $\vec{P}=\vec{\n}p$, where $p$ is a homogeneous
polynomial with an isolated singularity at the origin, sequence
$\vec{P}$ is regular, see \cite{AP}.
\end{itemize}

\subsection{Poisson cohomology in dimension 3}

Set $E:=\Ci(\R^3)$ and identify---as above---the spaces of
homogeneous multivector fields in $\R^3$, with the corresponding
component spaces: ${\cal X}^0(\R^3)\simeq{\cal X}^3(\R^3)\simeq E$
and ${\cal X}^1(\R^3)\simeq{\cal X}^2(\R^3)\simeq E^3$.

Let $\vec{\zL}=(\zL_1,\zL_2,\zL_3)\in E^3$ be a Poisson tensor and
$f\in E,\vec{X}\in E^3,\vec{B}\in E^3,T\in E$ a $0$-, $1$-, $2$-,
and $3$-cochain of the LP-complex. The following formul{\ae} for
the LP-coboundary operator $\p_{\vec{\zL}}$ can be obtained by
straightforward computations:
$$\begin{array}{lll}\p_{\vec{\zL}}^0f&=&\vec{\n}f\w\vec{\zL},\\
\p_{\vec{\zL}}^1\vec{X}&=&(\vec{\n}\cdot\vec{X})\vec{\zL}-\vec{\n}(\vec{X}\cdot\vec{\zL})+\vec{X}\w(\vec{\n}\w\vec{\zL}),\\
\p_{\vec{\zL}}^2\vec{B}&=&-(\vec{\n}\w\vec{B})\cdot\vec{\zL}-\vec{B}\cdot(\vec{\n}\w\vec{\zL}),\\
\p_{\vec{\zL}}^3T&=&0.\end{array}$$ If we denote the differential
detailed in Equation (\ref{RhamKoszul}) (resp. in Equation
(\ref{MultiplicKosul}) if $\vec{P}=\vec{\zL}$, in Equation
(\ref{MultiplicKosul}) if $\vec{P}=\vec{\n}\w\vec{\zL}$) by ${\cal
K}$ (resp. $\zk'$, $\zk''$), we get
\begin{equation}\p_{\vec{\zL}}^0=\zk'{\cal K},\p_{\vec{\zL}}^1=\zk'{\cal
K}-{\cal K}\zk'+\zk'',\p_{\vec{\zL}}^2=-\zk'{\cal
K}-\zk'',\p_{\vec{\zL}}^3=0.\label{LPCompKoszul}\end{equation}

As aforementioned, investigations are confined in this paper to
quadratic Poisson tensors and polynomial (or formal) LP-cochains.
If structure $\vec{\zL}$ is $K$-exact, i.e., in view of notations
due to the elimination of the module basis of multivector fields,
$\vec{\zL}=\vec{\n}p\;\;(p\in{\cal
S}^3\R^{3*})\Leftrightarrow\vec{\n}\w\vec{\zL}=0$, homology
operator $\zk''$ vanishes. If, moreover, $p$ has an isolated
singularity ({\small IS}), not only the $K^*$-complex associated
with ${\cal K}$ is exact up to injectivity of ${\cal
K}=\vec{\n}(.)$, but also the $K_*$-complex associated with $\zk'$
is, see above, acyclic up to surjectivity of
$\zk'=(.)\cdot\vec{\zL}$. In \cite{AP}, the author has computed
inter alia the LP-cohomology for a weight-homogeneous polynomial
$p$ with an {\small IS}.

Below, we describe a generic cohomological technique for {\small
SRMI} Poisson tensors in a finite-dimensional vector space. This
approach extends Formul{\ae} (\ref{LPCompKoszul}) to dimension $n$
and reduces simultaneously the LP-coboundary operator $\p_{\zL}$
to a single Koszul differential.

\subsection{Poisson cohomology in dimension n} \label{SecPoissCoho}

We denote by $L$ the matrix of maximal minors of a matrix
$\E\in\op{gl}(n,\R^{n*})$ (or of a matrix with entries in a field
$\F$ of non-zero characteristic), so $L_{ij}$ is the minor of $\E$
obtained by cancellation of line $i$ and column $j$. More
generally, if $\zn=\{1,\ldots,n\}$,
$\mathbf{i}=(i_1,\ldots,i_m)\in\zn^m$
($i_1<\ldots<i_m,\;m\in\{1,\ldots,n\}$), we denote by
$\mathbf{I}=(I_1,\ldots,I_{n-m})$ the complement of $\mathbf{i}$
in $\zn$. If $\mathbf{j}=(j_1,\ldots,j_m)$ is an $m$-tuple similar
to $\mathbf{i}$, we denote by $L_{\mathbf{i}\mathbf{j}}$ the minor
of $\E$ obtained by cancellation of the lines $\mathbf{i}$ and the
columns $\mathbf{j}$, and by $L^{\mathbf{i}\mathbf{j}}$ the minor
of $\E$ at the intersections of lines $\mathbf{i}$ and columns
$\mathbf{j}$. Hence,
$L_{\mathbf{i}\mathbf{j}}=L^{\mathbf{I}\mathbf{J}}$ and
$L_{\mathbf{I}\mathbf{J}}=L^{\mathbf{i}\mathbf{j}}.$ Moreover,
$D=\op{det}\E\in{\cal S}^n\R^{n*}$ is the determinant of $\E$,
${\cal L}$ stands for the matrix of maximal minors of
$L\in\op{gl}(n,{\cal S}^{n-1}\R^{n*})$, and we apply the just
introduced notations ${\cal L}_{\mathbf{i}\mathbf{j}}$ and ${\cal
L}^{\mathbf{i}\mathbf{j}}$ also to ${\cal L}$. Eventually, as
already mentioned above, $\mathbf{L}$ denotes the matrix of
algebraic maximal
minors of $\E.$\\

\noindent {\bf Remark}. In the following, we systematically assume
that $D\neq 0$, i.e. that polynomial $D$ does not vanish
everywhere.

\begin{lem}\label{LemMinors} For any $m\in\{1,\ldots,n-1\}$ and for any $\mathbf{i}=
(i_1,\ldots,i_m)$, $\mathbf{j}=(j_1,\ldots,j_m)$ as above, we have
$${\cal L}_{\mathbf{i}\mathbf{j}}=D^{n-m-1}L^{\mathbf{i}\mathbf{j}}\mbox{  and  }
{\cal L}^{\mathbf{i}\mathbf{j}}=D^{m-1}L_{\mathbf{i}\mathbf{j}}.$$
The first $($ resp. second $)$ equation also holds for $m=0$ $($
resp. $m=n$ $)$. In this case it just means that
$\op{det}L=D^{n-1}$.\end{lem}

{\it Proof}. Of course, the second statement is nothing but a
reformulation of the first. We prove the first assertion by
induction on $m$. For $m=n-1$, the assertion is obvious. Indeed,
the both sides coincide with the element $L_{IJ}$ of L at the
intersection of the line $I$ and the column $J$. Assume now that
the equation holds true for $2\le m\le n-1$ and take any
$\mathbf{i}=(i_1,\ldots,i_{m-1})$ and
$\mathbf{j}=(j_1,\ldots,j_{m-1})$ of length $m-1$. Let $i_m$ be an
(arbitrary) element of $(n-m+1)$-tuple $\mathbf{I}$. We will also
have to consider the $m$-tuple
$\mathbf{\underline{i}}=(i_1,\ldots,i_m,\ldots,i_{m-1})$, where
the elements have of course been written in the natural order
$i_1<\ldots<i_m<\ldots<i_{m-1}$. The rank of $i_m$ inside
$\mathbf{I}$ and $\mathbf{\underline{i}}$ will be denoted be
$r_{\mathbf{I}}(i_m)$ and $r_{\mathbf{\underline{i}}}(i_m)$
respectively. Using these notations, we get $${\cal
L}_{\mathbf{i}\mathbf{j}}={\cal
L}^{\mathbf{I}\mathbf{J}}=\sum_{j_m\in\mathbf{J}}(-1)^{r_{\mathbf{I}}(i_m)+r_{\mathbf{J}}(j_m)}L_{i_mj_m}{\cal
L}_{\mathbf{\underline{i}}\mathbf{\underline{j}}}.$$ Applying the
induction assumption, we see that ${\cal
L}_{\mathbf{\underline{i}}\mathbf{\underline{j}}}=D^{n-m-1}L^{\mathbf{\underline{i}}\mathbf{\underline{j}}}$,
so that $${\cal
L}_{\mathbf{i}\mathbf{j}}=D^{n-m-1}\sum_{j_m\in\zn}(-1)^{r_{\mathbf{I}}(i_m)+r_{\mathbf{J}}(j_m)}L_{i_mj_m}
\sum_{\zs\in{\cal
P}(\mathbf{\underline{j}})}\op{sign}\zs\,\E_{i_1\zs_{j_1}}\ldots\E_{i_m\zs_{j_m}}\ldots\E_{i_{m-1}\zs_{j_{m-1}}},$$
where ${\cal P}(\mathbf{\underline{j}})$ is the permutation group
of $\mathbf{\underline{j}}$, and where the first sum could be
extended to all $j_m\in\zn$, as for $j_m\in\mathbf{j}$ the last
determinant vanishes. It is clear that we obtain all the
permutations $\zs$ of $\mathbf{\underline{j}}$, if we assign $j_m$
to $i_p$ ($p\in\{1,\ldots,m\}$) and, for each choice of $p$, all
the permutations $\zm\in{\cal P}(\mathbf{j})$ to the remaining
subscripts $i_q$. Observe that the signature of the permutation
$\zs$ that associates $j_m$ with $i_p$ and permutes $\mathbf{j}$
by $\zm$, is
$\op{sign}\zs=(-1)^{r_{\mathbf{\underline{i}}}(i_p)-r_{\mathbf{\underline{j}}}(j_m)}\op{sign}\zm$.
Hence, we get $$\begin{array}{c}{\cal
L}_{\mathbf{i}\mathbf{j}}=\\D^{n-m-1}
\sum_{p=1}^m\sum_{\zm\in{\cal
P}(\mathbf{j})}\op{sign}\zm\,\E_{i_1\zm_{j_1}}\ldots\widehat{\E_{i_pj_m}}\ldots\\\E_{i_{m-1}\zm_{j_{m-1}}}
\sum_{j_m\in\zn}(-1)^{r_{\mathbf{I}}(i_m)+r_{\mathbf{\underline{i}}}(i_p)+r_{\mathbf{J}}(j_m)-r_{\mathbf{\underline{j}}}(j_m)}\E_{i_pj_m}L_{i_mj_m}.\end{array}$$
Remark now that the exponent of $-1$ can be replaced by
$$r_{\mathbf{\underline{i}}}(i_p)+r_{\mathbf{\mathbf{\underline{i}}}}(i_m)+r_{\mathbf{\mathbf{\underline{i}}}}(i_m)
+r_{\mathbf{I}}(i_m)+r_{\mathbf{\underline{j}}}(j_m)+r_{\mathbf{J}}(j_m)\sim
r_{\mathbf{\underline{i}}}(i_p)+r_{\mathbf{\mathbf{\underline{i}}}}(i_m)+i_m+j_m.$$
Thus, the last sum reads
$(-1)^{r_{\mathbf{\underline{i}}}(i_p)+r_{\mathbf{\mathbf{\underline{i}}}}(i_m)}\sum_{j_m\in\zn}
(-1)^{i_m+j_m}\E_{i_pj_m}L_{i_mj_m}$. If $p\neq m$, this sum
vanishes, and if $p=m$ it coincides with determinant $D$.
Eventually, we find $${\cal
L}_{\mathbf{i}\mathbf{j}}=D^{n-m}\sum_{\zm\in{\cal
P}(\mathbf{j})}\op{sign}\zm\,\E_{i_1\zm_{j_1}}\ldots\widehat{\E_{i_mj_m}}\ldots\E_{i_{m-1}\zm_{j_{m-1}}}=D^{n-m}L^{\mathbf{i}\mathbf{j}}.\quad\rule{1.5mm}{2.5mm}$$

\begin{defi} Let $Y_i=\sum_r\E_{ir}\p_r$ be $n$ linear vector fields in $\R^n$. Set
$${\cal R}=\oplus_{p=0}^n{\cal
R}^p=\oplus_{p=0}^n\R[[x_1,\ldots,x_n]]\otimes\w^p_n\la\vec{\p}\ra$$
and $${\cal P}=\oplus_{p=0}^n{\cal
P}^p=D^{-1}\oplus_{p=0}^n\R[[x_1,\ldots,x_n]]\otimes\w^p_n\la\vec{Y}\ra,$$
where $D=\op{det}\E$ and where $\w_n^p\la\vec{\p}\ra$ and
$\w_n^p\la\vec{Y}\ra$ are the terms of degree $p$ of the Grassmann
algebras on generators $\vec{\p}=(\p_1,\ldots,\p_n)$ and
$\vec{Y}=(Y_1,\ldots,Y_n)$ respectively. Space ${\cal R}$
$(\mbox{resp. }{\cal P})$ is the space of {\sl real} formal
LP-cochains (resp. {\sl potential} formal LP-cochains).
\end{defi}

\noindent {\bf Remark}. The space of polyvector fields $Y_{\mathbf
k}=Y_{k_1\ldots k_p}=Y_{k_1}\w\ldots\w Y_{k_p}$ $(k_1<\ldots<k_p,
p\in\{0,\ldots,n\})$ with coefficients in the quotients by $D$ of
formal power series in $(x_1,\ldots,x_n)$, is a concrete model of
space ${\cal P}$. Indeed, observe first that these spaces are
bigraded by the ``exterior degree'' $p$ and the (total)
``polynomial degree'', say $r$. If such a polyvector field
vanishes, its homogeneous terms
$D^{-1}\sum_{\mathbf{k}}P^{\mathbf{k}r}Y_{\mathbf{k}}$
$(P^{\mathbf{k}r}\in{\cal S}^r\R^{n*})$ vanish. If we decompose
the $Y_i$ $(i\in\{1,\ldots,n\})$ in the natural basis $\p_i$, we
immediately see that the sums
$\sum_{\mathbf{k}}L^{\mathbf{k}\mathbf{i}}P^{\mathbf{k}r}$ vanish
for all $\mathbf{i}=(i_1,\ldots,i_p)$ $(i_1<\ldots<i_p).$ Since
these sums can be viewed as the product of a matrix with
polynomial entries and the column made up by the $P^{\cdot r}$,
the column vanishes outside the vanishing set $V$ of the
homogeneous polynomial determinant of this matrix. As the
complement of (the conic closed) subset $V$ of $\R^n$ is dense in
$\R^n$, the polynomials $P^{\mathbf{k}r}$ vanish everywhere.

\begin{theo}\label{injection} (i) There is a canonical non surjective injection $i:{\cal
R}\raa{\cal P}$ from ${\cal R}$ into ${\cal P}$.\\
(ii) A homogeneous potential cochain
$D^{-1}\sum_{\mathbf{k}}P^{\mathbf{k}r}Y_{\mathbf{k}}$ $[$of
bidegree $(p,r)$$]$ is real if and only if the $[$$n!/p!(n-p)!$$]$
homogeneous polynomials
$\sum_{\mathbf{k}}L^{\mathbf{k}\mathbf{i}}P^{\mathbf{k}r}$ $[$of
degree $p+r$$]$ are divisible by $D$ $($for $p=0$ this condition
means that $P^r$ be divisible by $D$$)$.
\end{theo}

{\it Proof}. Take a real cochain
$C^p=\sum_{\mathbf{i}}\varsigma^{\mathbf{i}}\p_{\mathbf{i}}\in{\cal
R}^p$, where, as above, $\mathbf{i}=(i_1,\ldots,i_p)$,
$i_1<\ldots<i_p$. As $\p_j=D^{-1}\sum_k\mathbf{L}_{kj}Y_k$, we get
$$\p_{\mathbf{i}}=D^{-p}\sum_{k_1,\ldots,k_p}\mathbf{L}_{k_1i_1}\ldots\mathbf{L}_{k_pi_p}Y_{k_1\ldots
k_p}=D^{-p}\sum_{k_1<\ldots<k_p}\lp\sum_{\zs\in{\cal
P}(\mathbf{k})}\op{sign}\zs\,\mathbf{L}_{\zs_{k_1}i_1}\ldots\mathbf{L}_{\zs_{k_p}i_p}\rp
Y_{k_1\ldots k_p}.$$ If $\mid\!\mathbf{i}\!\mid=\sum_{j=1}^p i_j$,
it follows from Lemma \ref{LemMinors}, that the determinant in the
above bracket is given by
$$(-1)^{\mid\mathbf{i}\mid+\mid\mathbf{k}\mid}{\cal L}^{\mathbf{k}\mathbf{i}}=(-1)^{\mid\mathbf{i}\mid+\mid\mathbf{k}\mid}D^{p-1}L_{\mathbf{k}\mathbf{i}},$$
so that
$$C^p=D^{-1}\sum_{\mathbf{k}}\lp\sum_{\mathbf{i}}(-1)^{\mid\mathbf{i}\mid+\mid\mathbf{k}\mid}L_{\mathbf{k}\mathbf{i}}\,\zvs^{\mathbf{i}}\rp
Y_{\mathbf{k}},$$ where the {\small RHS} is in ${\cal P}^p.$

Point (ii) is a direct consequence of the preceding remark.
\rule{1.5mm}{2.5mm}\\

\noindent {\bf Remark}. In view of this theorem, the bigrading
${\cal P}=\oplus_{p=0}^n\oplus_{r=0}^{\infty}{\cal P}^{pr}$,
defined on ${\cal P}$ by the exterior degree and the polynomial
degree, induces a bigrading ${\cal
R}=\oplus_{p=0}^n\oplus_{r=0}^{\infty}{\cal R}^{pr}$ on ${\cal
R}$.\\

Consider now a quadratic Poisson tensor $\zL$ in $\R^n$. In the
following, we {\it assume} that $\zL$ is {\small SRMI}, and more
precisely that there are $n$ mutually commuting linear vector
fields $Y_i=\sum_{r=1}^n\E_{ir}\p_r$, $\,\E\in\op{gl}(n,\R^{n*})$,
such that $D=\op{det}\,\E\neq0$ and
$$\zL=\sum_{i<j}\za^{ij}Y_{ij}\;\;\;(\za^{ij}\in\R).$$

\begin{prop}\label{JointEigen} The determinant $D=\op{det}\E\in{\cal
S}^n\R^{n*}\setminus\{0\}$ of $\,\E$ is the unique joint
eigenvector of the $Y_i$ with eigenvalues $\op{div}Y_i\in\R$,
i.e., $D$ is, up to multiplication by nonzero constants, the
unique nonzero polynomial of $\R^n$ that verifies
$$Y_iD=(\op{div}Y_i)D,\forall i\in\{1,\ldots,n\}.$$ Moreover, if
$D=D_1D_2$, where $D_1\in{\cal S}^{n_1}\R^{n*}$ and $D_2\in{\cal
S}^{n_2}\R^{n*}$ $(n_1+n_2=n)$ are two polynomials without common
divisor, these factors $D_1$ and $D_2$ are also joint eigenvectors.
If $\zl_i$ and $\zm_i$ denote their eigenvalues, we have
$\zl_i+\zm_i=\op{div}Y_i$.
\end{prop}

{\it Proof}. Set $Y_i=\sum_r\E_{ir}\p_r=\sum_{rs}a_{ir}^sx_s\p_r$,
$a_{ir}^s\in\R$. Note first that
$Y_{i}(\ell_{jr})=\sum_{t}a_{jr}^{t}\ell_{it}$, and that
$[Y_{i},Y_{j}]=0$ means $Y_{i}(\ell _{jr})=Y_{j}(\ell _{ir})$, for
all $i,j,r\in \{1,\ldots ,n\}$. If ${\cal P}_n$ denotes the
permutation group of $\{1,\ldots,n\}$, we then get
\begin{eqnarray*} Y_{i}D &=&\sum_{k=1}^{n}\sum_{\sigma \in
\mathcal{P}_n}\op{sign}\!\sigma\;\ell
_{\sigma _{1}1}\ldots Y_{i}(\ell _{\sigma _{k}k})\ldots\ell _{\sigma _{n}n} \\
&=&\sum_{k=1}^{n}\sum_{\sigma \in
\mathcal{P}_n}\op{sign}\!\sigma\;\ell _{\sigma
_{1}1}\ldots Y_{\sigma _{k}}(\ell _{ik})\ldots\ell _{\sigma _{n}n} \\
&=&\sum_{k,t=1}^{n}a_{ik}^{t} \sum_{\sigma \in
\mathcal{P}_n}\op{sign}\!\sigma\; \ell _{\sigma _{1}1}\ldots\ell
_{\sigma _{k}t}\ldots\ell _{\sigma _{n}n}.
\end{eqnarray*} This last sum vanishes if $k\neq t$ since two columns coincide in this determinant.
Eventually, we have $$Y_iD=\left(\sum_k a_{ik}^{k}\right)D=\left(\op{div}Y_i\right)D.$$ \\
\indent As for uniqueness, suppose that there is another
polynomial $P\in{\cal S}\R^{n*}\setminus\{0\}$, such that
$Y_iP=(\op{div}Y_i)P$, for all $i\in\{1,\ldots,n\}$. Then
$Y_i\left(P/D\right)=0$ in $Z=\{x\in\R^n,D(x)\neq 0\}$ and the
same reasoning as in the proof of Theorem \ref{1} allows
concluding that there exists $\alpha\in\R^*$ such that $P=\alpha
D$.

The assertion concerning the factorization of $D$ is easily
understood. Indeed, since
$((\op{div}Y_i)D_1-Y_iD_1)D_2=D_1(Y_iD_2)$ and as the polynomials
$D_1$ and $D_2$ have no common divisor, $Y_iD_2=PD_2$ and
$(\op{div}Y_i)D_1-Y_iD_1=QD_1$, where $P=Q$ is a polynomial.
Looking at degrees, we immediately see that $P=Q$ is necessarily
constant.
\rule{1.5mm}{2.5mm}\\

\noindent {\bf Remark}. Observe that the eigenvalues
$\op{div}Y_i$, $i\in\{1,\ldots,n\}$, cannot vanish simultaneously.
Indeed, in this case, polynomial $D\in{\cal
S}^n\R^{n*}\setminus\{0\},$ $n\in\N^*$, vanishes everywhere.

\begin{defi} The complex $$0\raa {\cal R}^0\raa {\cal R}^1\raa\ldots\raa {\cal R}^n\raa 0$$
with differential $\p_{\zL}=[\zL,.]_{\op{SN}}$, is the formal
LP-complex of Poisson tensor $\zL\in{\cal
S}^2\R^{n*}\otimes\w^2\R^n$. We denote the corresponding
cohomology groups by $LH^*({\cal R},\zL)$.\end{defi}

The next theorem shows that if the cochains $C\in{\cal R}$ are
read as $C=iC\in{\cal P}$, the LP-differential assumes a
simplified shape.

\begin{theo}\label{PoissKoszCob} Set $\zL=\sum_{i<j}\za^{ij}Y_{ij}$, $\za^{ji}=-\za^{ij}$, and $X_i=\sum_{j\neq
i}\za^{ij}Y_j$.\\
(i) Let
$$C=D^{-1}\sum_{\mathbf{k}}P^{\mathbf{k}r}Y_{\mathbf{k}}\in{\cal P}^{pr}$$
be a homogeneous potential cochain. The LP-coboundary of $C$ is
given by \be\p_{\zL}C=\sum_{\mathbf{k}i}X_i\!\lp
D^{-1}P^{\mathbf{k}r}\rp Y_i\w Y_{\mathbf{k}}
=D^{-1}\sum_{\mathbf{k}i}\lp
X_i-\zd_i\op{id}\rp\!(P^{\mathbf{k}r})\,Y_i\w
Y_{\mathbf{k}}\in{\cal P}^{p+1,r},\label{SRMICobFin}\ee
where $\zd_i=\op{div}X_i\in\R$.\\

\noindent (ii) The LP-coboundary operator $\p_{\zL}$ endows ${\cal
P}$ with a differential complex structure, and preserves the
polynomial degree $r$. This LP-complex of $\zL$ over ${\cal P}$
contains the LP-complex $({\cal R},\p_{\zL})$ of $\zL$ over ${\cal
R}$ as a differential sub-complex.
\end{theo}

{\it Proof}. Note first that if $C=f \mathbf{Y},$ where $f$ a
function and $\mathbf{Y}$ a wedge product of vector fields $Y_k$, we
get
\be\p_{\zL}(f\mathbf{Y})=[\zL,f\mathbf{Y}]_{\op{SN}}=[\zL,f]_{\op{SN}}\w\mathbf{Y},\label{SRMICob1}\ee
since the $Y_k$ are mutually commuting. However,
\be[\zL,f]_{\op{SN}}=\sum_{i<j}\za^{ij}\lp
(Y_jf)Y_i-(Y_if)Y_j\rp=\sum_i(\sum_{j\neq i}\za^{ij}\,Y_jf)
Y_i=\sum_i(X_if)Y_i.\label{SRMICob2}\ee When combining Equations
(\ref{SRMICob1}) and (\ref{SRMICob2}), we get the first part of
Equation (\ref{SRMICobFin}), whereas its second part is the
consequence of Proposition \ref{JointEigen}.  \rule{1.5mm}{2.5mm}

\begin{cor} The LP-cohomology groups of $\zL$ over ${\cal R}$ and ${\cal P}$ are
bigraded, i.e.
$$LH({\cal R},\zL)=\oplus_{r=0}^{\infty}\oplus_{p=0}^nLH^{pr}({\cal R},\zL)\quad\mbox{and}\quad
LH({\cal P},\zL)=\oplus_{r=0}^{\infty}\oplus_{p=0}^nLH^{pr}({\cal
P},\zL),$$ where for instance $LH^{pr}({\cal P},\zL)$ is defined
by
$$LH^{pr}({\cal P},\zL)=\op{ker}(\p_{\zL}:{\cal P}^{pr}\raa {\cal
P}^{p+1,r})/\op{im}(\p_{\zL}:{\cal P}^{p-1,r}\raa{\cal P}^{pr}).$$
\end{cor}

In the following we deal with the terms $LP^{*r}({\cal
P},\zL)=\oplus_{p=0}^nLP^{pr}({\cal P},\zL)$ of the LP-cohomology
over ${\cal P}$ and with the corresponding part of LP-cohomology
the subcomplex ${\cal R}$.

\begin{theo}\label{PoissKosz} Let $E_r$ be the real finite-dimensional vector space ${\cal
S}^r\R^{n*},$ and let
$\vec{X}_{\delta}:=(X_1-\delta_1\op{id},\ldots,X_n-\delta_n\op{id})$,
$\zd_i=\op{div}X_i$, be the $n$-tuple of the commuting linear
operators $X_i-\delta_i\op{id}$ on $E_r$ defined in Theorem
\ref{PoissKoszCob}. The LP-cohomology space $LH^{*r}({\cal
P},\zL)$ of $\zL$ over ${\cal P}$ coincides with the Koszul
cohomology space $KH^*(\vec{X}_{\zd},E_r)$ associated with
$\vec{X}_{\zd}$ on $E_r$: $$LH^{*r}({\cal P},\zL)\simeq
KH^*(\vec{X}_{\zd},E_r).$$
\end{theo}

{\it Proof}. Direct consequence of result $\p_{\zL}=\sum_i
(X_i-\delta_i\op{id})\otimes e_{Y_i}$ proved in Theorem \ref{PoissKoszCob}.
\rule{1.5mm}{2.5mm}\\

In order to study the LP-cohomology group $LH^{.r}({\cal R},\zL)$
of the quadratic Poisson tensor $\zL$ over the formal cochain
space ${\cal R}$, we introduce a long cohomology exact
sequence.\\

Let ${\cal S}^{pr}$ be a complementary vector subspace of ${\cal
R}^{pr}$ in ${\cal P}^{pr}$: ${\cal P}^{pr}={\cal
R}^{pr}\oplus{\cal S}^{pr}.$ Space ${\cal
S}=\oplus_{r=0}^{\infty}\oplus_{p=0}^n{\cal S}^{pr}$ can easily be
promoted into the category of differential spaces. Indeed, denote
by $p_{\cal R}$ and $p_{\cal S}$ the projections of ${\cal P}$
onto ${\cal R}$ and ${\cal S}$ respectively, and set for any
$s\in{\cal S}$,
\[\zf s=p_{\cal R}\p_{\zL}s,\tilde{\p}_{\zL}s=p_{\cal S}\p_{\zL}s.\]

\begin{prop}
(i) The endomorphism $\tilde{\p}_{\zL}\in \op{End}_{\R}{\cal S}$
is a differential on ${\cal S}$, which has weight $(1,0)$ with
respect to the bigrading of ${\cal S}$, i.e.
$\tilde{\p}_{\zL}:{\cal S}^{pr}\raa {\cal S}^{p+1,r}$.

\noindent (ii) The linear map $\zf\in \op{Hom}_{\R}({\cal S},{\cal
R})$ is an anti-homomorphism of differential spaces from $({\cal
S},\tilde{\p}_{\zL})$ into $({\cal R},\p_{\zL})$. Its weight with
respect to the bidegree is $(1,0)$, i.e. $\zf:{\cal S}^{pr}\raa
{\cal R}^{p+1,r}$.

\noindent (iii) The sequence $0\raa {\cal
R}\stackrel{i}{\raa}{\cal P}\stackrel{p_{\cal S}}{\raa}{\cal
S}\raa 0$ is a short exact sequence of homomorphisms of
differential spaces, which preserve the bidegree. It induces an
exact triangle in cohomology, whose connecting homomorphism
$\zf_{\sharp}$ is canonically implemented by $\zf$. If
$LH^{pr}({\cal S},\tilde{\zL})$ denotes the degree $(p,r)$ term of
the cohomology space of the complex $({\cal S},\tilde{\p}_{\zL})$,
we have $\zf_{\sharp}:LH^{pr}({\cal S},\tilde{\zL})\raa
LH^{p+1,r}({\cal R},\zL)$.

\noindent (iv) The sequence
\[\begin{array}{l}0\raa LH^{0r}({\cal
R},\zL)\stackrel{i_{\sharp}}{\raa}\ldots\\\;\;\;\;\;\;\;
\stackrel{\zf_{\sharp}}{\raa} LH^{pr}({\cal
R},\zL)\stackrel{i_{\sharp}}{\raa}LH^{pr}({\cal P},\zL)
\stackrel{(p_{\cal S})_{\sharp}}{\raa}LH^{pr}({\cal
S},\tilde{\zL})\stackrel{\zf_{\sharp}}{\raa}LH^{p+1,r}({\cal
R},\zL)\stackrel{i_{\sharp}}{\raa}\ldots\stackrel{(p_{\cal
S})_{\sharp}}{\raa}LH^{nr}({\cal S},\tilde{\zL})\raa
0\end{array}\] is a long exact cohomology sequence of vector space
homomorphisms.

\noindent (v) If $\op{ker}^{pr}\zf_{\sharp}$ and
$\op{im}^{p+1,r}\zf_{\sharp}$ denote the kernel and the image of
the restricted map $\zf_{\sharp}:LH^{pr}({\cal S},\tilde{\zL})\raa
LH^{p+1,r}({\cal R},\zL)$, we have\be LH^{pr}({\cal R},\zL) \simeq
LH^{p-1,r}({\cal
S},\tilde{\zL})/\op{ker}^{p-1,r}\zf_{\sharp}\oplus LH^{pr}({\cal P
},\zL)/\op{ker}^{pr}\zf_{\sharp}.\label{PoissRealPoten}\ee
\end{prop}

{\it Proof.} Statements (i) and (ii) are direct consequences of
equation $\p_{\zL}^2=0.$ For (iii), we only need check that linear
map $\zf_{\sharp}$ coincides with the connecting homomorphism,
what is obvious. Eventually, assertion (v) is a corollary of
exactness of the long cohomology sequence. \rule{1.5mm}{2.5mm}\\

We now identify the $\cal{S}$-cohomology with a relative
cohomology. Several concepts of relative cohomology can be met in
literature. Below, we use the following definition.

\begin{defi} Let $V$ be a vector space endowed with a differential
$\p$, and let $W$ be a $\p$-closed subspace of $V$. Denote by
$\overline{\p}$ the differential canonically induced by $\p$ on
the quotient space $V\slash W$. The cohomology of the differential
space $(V\slash W,\overline{\p})$ is called the relative
cohomology of $(V,W,\p)$. It is denoted by $H(V,W,\p)$.
\end{defi}

\begin{prop}\label{SuppRelat} The cohomology induced by $\p_{\zL}$ on ${\cal S}$ $($i.e. the
cohomology of differential space $({\cal S},\tilde{\p}_{\zL})$$)$
coincides with the relative cohomology of $({\cal P},{\cal
R},\zL)$ $($i.e. the cohomology of space $({\cal P}\slash{\cal
R},\overline{\p}_{\zL})$$)$:
$$LH({\cal S},\tilde{\zL})\simeq LH({\cal P},{\cal R},\zL).$$\end{prop}

{\it Proof}. It suffices to note that the vector space isomorphism
$\psi:{\cal P}\slash {\cal R}\ni [\zp]\raa p_{\cal S}\zp\in{\cal
S}$ intertwines the differentials $\overline{\p}_{\zL}$ on ${\cal
P}\slash{\cal R}$ and $\tilde{\p}_{\zL}$ on ${\cal S}$. \rule{1.5mm}{2.5mm}\\

\noindent {\bf Remark}. In view of this proposition it is clear
that ${\cal S}$-cohomology is independent of the chosen splitting
${\cal P}={\cal R}\oplus{\cal S}$.

\begin{theo}\label{PoissKoszRel} The LP-cohomology groups of a {\small SRMI} Poisson tensor $\zL,$
over the space ${\cal R}$ of cochains with coefficients in the
formal power series, are given by \[LH^{pr}({\cal R},\zL) \simeq
LH^{pr}({\cal P },\zL)/\op{ker}^{pr}\zf_{\sharp}\oplus
LH^{p-1,r}({\cal P},{\cal R},\zL)/\op{ker}^{p-1,r}\zf_{\sharp},\]
where the above-introduced notations have been used. \end{theo}

{\it Proof}. Reformulation of Equation (\ref{PoissRealPoten}) and
Proposition \ref{SuppRelat}. \rule{1.5mm}{2.5mm}\\

\noindent {\bf Remark}.{\it This theorem reduces computation of
the formal LP-cohomology groups $LH^{pr}({\cal R},\zL)$, basically
to the Koszul cohomology groups $LH^{pr}({\cal P},\zL)\simeq
KH^p(\vec{X}_{\zd},E_r)$ associated to the
afore-detailed operators $\vec{X}_{\zd}$ on
$E_r={\cal S}^r\R^{n*}$ induced by the considered {\small SRMI} tensor, and to the relative cohomology groups
$LH^{p-1,r}({\cal P},{\cal R},\zL)$. It thus highlights the link
between Poisson and Koszul cohomology. Let us mention that we
showed in \cite{MP}, via explicit computations in $\R^3$, that
${\cal P}$-cohomology $($now identified as Koszul cohomology$)$
and ${\cal S}$-cohomology $($or relative cohomology$)$ are less
intricate than Poisson cohomology.}\\

The remark concerning the comparative simplicity of the ${\cal
P}$-cohomology can be easily understood.

Observe that any {\small SRMI} Poisson tensor
$\zL=\sum_{i<j}\za^{ij}Y_{ij},$ $\za^{ij}\in\R,$ with
$Y_i=\sum_r\E_{ir}\p_r$ and $\E_{ir}\in\R^{n*}$, reads, locally in
$\{D:=\op{det}\E\neq 0\}\subset\{\zL\neq 0\}\subset\R^n$,
$$\zL=\sum_{i<j}\za^{ij}\p_{s_is_j}\;\;(\za^{ij}\in\R),$$ where $(s_1,\ldots,s_n)$
are local coordinates. As the $Y_i$ mutually commute, the
statement is a direct consequence of the ``straightening theorem
for vector fields''. For instance, for structure
$\zL=2a\,Y_{23}+b\,Y_{31}+a\,Y_{12},$ where
$Y_1=x\p_1+y\p_2,\,Y_2=x\p_2-y\p_1,\,Y_3=z\p_3,$ and
$D=(x^2+y^2)z$, see Theorem \ref{ClassTheo}, the local
(non-polynomial) coordinate transformation
$$x=e^s\op{cos}\theta,y=e^s\op{sin}\theta,z=-e^{-t}$$ leads to
$Y_1=\p_s,Y_2=\p_{\theta},Y_3=\p_{t},$ and $\zL=2a\p_{\theta
t}+b\p_{ts}+a\p_{s\theta}.$

Hence, {\it locally} in a dense open subset of $\R^n$, there are
coordinate systems or bases in which tensor $\zL$ has {\it
constant coefficients}. The ${\cal P}$-cohomology $LH^{*r}({\cal
P},\zL)$ however, is the LP-cohomology in the extended space
${\cal P}^{*r}=D^{-1}\oplus_{p=0}^n{\cal
S}^r\R^{n*}\otimes\w^p_n\la\vec{Y}\ra$, which admits the {\it
global} basis $\vec{Y}=(Y_1,\ldots,Y_n)$ in which structure $\zL$
has {\it constant coefficients}. This is what makes ${\cal
P}$-cohomology particularly convenient.

\section{Koszul cohomology in a finite-dimensional vector space}

In view of the above remark regarding the basic ingredients of
LP-cohomology of {\small SRMI} tensors of $\R^n$, we take in this
section an interest in the Koszul cohomology space
$KH^*(\vec{X}_{\zl},E)$ associated to operators
$\vec{X}_{\zl}:=(X_1-\zl_1\op{id},\ldots,X_n-\zl_n\op{id})$ made
up of commuting linear transformations $\vec{X}:=(X_1,\ldots,X_n)$
of a finite-dimensional real vector space $E$ and a point
$\vec{\zl}:=(\zl_1,\ldots,\zl_n)\in\R^n$. However, Koszul
cohomology is known to be closely connected with Spectral
Theory---a fundamental principle of multivariate operator theory
is that all essential spectral properties of operators $\vec{X}$
in a complex space should be understood in terms of properties of
the Koszul complex induced by $\vec{X}_{\zl}$,
$\vec{\zl}\in\C^n$---so that the natural framework for
investigations on Koszul cohomology is the complex setting.

\begin{prop}\label{ComplexCoho} Let $(E,\p)$ be a differential space over $\R$, and denote
by $(E^{\C},\p^{\C})$ its complexification. The complexification
$H^{\C}(E,\p)$ of the cohomology space of $(E,\p)$ and the
cohomology $H(E^{\C},\p^{\C})$ of differential space
$(E^{\C},\p^{\C})$, are canonically isomorphic:
$$H(E^{\C},\p^{\C})\simeq H^{\C}(E,\p).$$
\end{prop}

{\it Proof}. Obvious. \rule{1.5mm}{2.5mm} 

\begin{prop}\label{ComplexCohoKosz} If $\vec{X}\in\op{End}_{\R}(E)$ are commuting $\R$-linear
transformations of a real vector space $E$, and if
$\vec{X}^{\C}\in\op{End}_{\C}(E^{\C})$ are the commuting
corresponding complexified $\C$-linear transformations of the
complexification $E^{\C}$ of $E$, the following isomorphism of
complex vector spaces holds:
$$KH^*(\vec{X}^{\C},E^{\C})\simeq KH^{*\C}(\vec{X},E).$$
\end{prop}

{\it Proof}. In view of Proposition \ref{ComplexCoho}, it suffices
to check that the complex $K^*(\vec{X}^{\C},E^{\C})$ is
effectively the complexification of the complex $K^*(\vec{X},E)$.
\rule{1.5mm}{2.5mm}\\

This proposition allows deducing our subject for investigation,
the Koszul cohomology $KH^*(\vec{X}_{\zl},E)$ (where $\vec{\zl}$
is a point of $\R^n$ and where $\vec{X}$ is an $n$-tuple of
commuting $\R$-linear operators of a finite-dimensional vector
space $E$ over $\R$), from its more natural counterpart over the
field of complex numbers.\\

Below, we use the concept of joint spectrum $\zs(\vec{X})$ of
commuting bounded linear operators $\vec{X}=(X_1,\ldots,X_n)$ on a
complex vector space $E$. There are a number of definitions of
such spectra in the literature; the considered spaces $E$ are
normed spaces, Banach spaces, or Hilbert spaces. Here we
investigate Koszul cohomology in finite dimension and need the
following characterizations of the elements of the joint spectrum
$\zs(\vec{X})$ (for a proof, we refer the reader to \cite{BR}):

\begin{prop}\label{CharFinDimJointSpec} Let $\vec{X}=(X_1,\ldots,X_n)$ be an n-tuple of commuting operators on a finite-dimensional
complex vector space $E$. Then the following statements are
equivalent for any fixed
$\vec{\zl}=(\zl_1,\ldots,\zl_n)\in\C^{n}$:
\begin{description}
  \item[(a)] $\vec{\zl}\in\sigma(\vec X)$
  \item[(b)] There exists a basis in $E$ with respect to which the matrices representing the $X_j$
  are all upper-triangular, and there exists an index $q$ $(1\leq q \leq \op{dim}E)$, such that $\zl_j$ is the $(q,q)$ entry of the matrix
  representing $X_j$, for $j\in\{1,\ldots,n\}$
  \item[(c)] There exists an index $q$ as in Item $(b)$ for every basis in $E$ with respect to which the matrices representing the $X_j$ are all upper-triangular
  \item[(d)] There exists a nonzero vector $x$ such that $X_jx=\zl_jx$, $\forall j\in\{1,\ldots,n\}$
   \item[(e)] There do not exist $Y_j$ in the subalgebra of $\op{End}_{\C}(E)$ generated by $\op{id}$ and $\vec X$,
   such that $$\sum_{j=1}^nY_j(X_j-\zl_j\op{id})=\op{id}$$
\end{description}
\end{prop}

In the following, we supply some results regarding Koszul
cohomology spaces. We use the same notations as above.

\begin{prop}\label{HomotGrass} Let $\w=\w_n\la \vec{\zh}\ra$ be the exterior algebra
on $n$ generators $\vec{\zh}$ over a field $\F$ of characteristic
$0$, and let $\vec{h}$ be dual generators, i.e.
$i_{h_k}\zh_{\E}=\p_{k\E}$. We then have the homotopy formula
$$e_{\zh_{\E}}i_{h_k}+i_{h_k}e_{\zh_{\E}}=\zd_{k\E}\op{id},$$ where
$i_{h_k}$ and $e_{\zh_{\E}}$ are the creation and annihilation
operators, respectively.
\end{prop}

{\it Proof}. Obvious. \rule{1.5mm}{2.5mm}

\begin{prop}\label{HomotKosz} Let $\vec{\cal X}\in\op{End}^{\times n}_{\F}(E)$ $($resp. $\vec{Y}\in\op{End}^{\times
n}_{\F}(E)$ $)$ be $n$ commuting linear operators $\vec{\cal X}$
$($resp. $\vec{Y}$$)$ on a vector space $E$ over $\F$. We denote
by ${\cal K}=\sum_{\E}{\cal X}_{\E}\otimes e_{\zh_{\E}}$ $($resp.
$\zk=\sum_kY_k\otimes i_{h_k}$$)$ the corresponding Koszul
cohomology $($resp. homology$)$ operator. The following
homotopy-type result holds: $${\cal K}\zk+\zk{\cal
K}=\lp\sum_{\E}Y_{\E}{\cal
X}_{\E}\rp\otimes\op{id}+\sum_{k\E}[{\cal X}_{\E},Y_{k}]\otimes
e_{\zh_{\E}}i_{h_{k}}.$$
\end{prop}

{\it Proof}. Direct consequence of Proposition \ref{HomotGrass}.
\rule{1.5mm}{2.5mm}

\begin{prop}\label{KoszCohoBasic} Let $\vec{X}\in\op{End}^{\times
n}_{\C}(E)$ be $n$ commuting endomorphisms of a finite-dimensional
complex vector space $E$, and let $\vec{\zl}\in\C^n$. Consider a
splitting $$E=E_1\oplus E_2$$ and denote by $i_{j}:E_j\raa E$
$($resp. $p_{j}:E\raa E_j$ $)$ the injection of $E_j$ into $E$
$($resp. the projection of $E$ onto $E_j$ $)$.

If $E_1$ is stable under the action of the operators $X_{\E}$,
i.e. $p_2\,X_{\E}\,i_1=0,$ and if $\vec{\zl}$ is not in the joint
spectrum $\zs(\vec{X}_{2})$ of the commuting operators $X_{\E
2}=p_2\,X_{\E}\,i_2\in\op{End}_{\C}(E_2)$, then any cocycle $C\in
E\otimes \w$ of the Koszul complex $K^*(\vec{X}_{\zl},E),$ where
$\vec{X}_{\zl}=\vec{X}-\vec{\zl}\op{id}_E$, is cohomologous to a
cocycle $C_1\in E_1\otimes\w,$ with $\w=\w_n\la\vec{\zh}\ra$.
\end{prop}

{\it Proof}. Observe first that if
$q(\vec{X})\in\C[X_1,\ldots,X_n]\subset\op{End}_{\C}(E)$ denotes a
polynomial in the $X_{\E}$, the compound map
$q(\vec{X})_{2}=p_2q(\vec{X})i_2$ coincides with the (same)
polynomial $q(\vec{X}_{2})\in\op{End}_{\C}(E_2)$ in the $X_{\E 2}$
($\star$). Indeed, due to stability of $E_1$, we have
$$p_2\,X_{\E}X_{k}\,i_2=p_2\,X_{\E}\,i_1\,p_1
\,X_{k}\,i_2+p_2\,X_{\E}\,i_2\,p_2\,X_{k}\,i_2=X_{\E 2}X_{k 2}.$$
This entails in particular that the $X_{\E 2}$ commute.

As $\vec{\zl}\notin\zs(\vec{X}_{2})$, Item (e) in Proposition
\ref{CharFinDimJointSpec} implies that there are $n$ operators
$\vec{Y}_{2}$ in the subalgebra of $\op{End}_{\C}(E_2)$ generated
by $\op{id_{E_2}}$ and $\vec{X}_{2}$, such that
\be\label{LeftInv2}\sum_{\E}Y_{\E 2}(X_{\E
2}-\zl_{\E}\op{id}_{E_2})=\op{id}_{E_2}.\ee Hence, for any $\E$,
$Y_{\E 2}=Q_{\E}(\vec{X}_2)$ is a polynomial in the $X_{k 2}$. Set
now $Y_{\E}=Q_{\E}(\vec{X})\in\op{End}_{\C}(E)$.

If applied to operators $\vec{X}_{\zl}$ and $\vec{Y}$, Proposition
\ref{HomotKosz} implies that
$$\lp\sum_{\E}Y_{\E}(X_{\E}-\zl_{\E}\op{id}_E)\rp\otimes\op{id}_{\w}+\sum_{k\E}[X_{\E}-\zl_{\E}\op{id}_E
,Y_{k}]\otimes e_{\zh_{\E}}i_{h_{k}}={\cal K}\zk+\zk{\cal K},$$
where ${\cal K}$ (resp. $\zk$) is the Koszul cohomology (resp.
homology) operator associated with $\vec{X}_{\zl}$ (resp.
$\vec{Y}$) on $E$. As $Y_{k}$ is a polynomial in the commuting
endomorphisms $X_{\E}$, the second term on the {\small LHS} of the
preceding equation vanishes. Hence, when evaluating both sides on
a cocycle $C=e\otimes w$ of cochain complex
$K^*(\vec{X}_{\zl},E)$, we get
$$\lp Q(\vec{X})(e)\rp w={\cal K}\zk(e\otimes w),$$ where
$Q(\vec{X})=\sum_{\E}Y_{\E}(X_{\E}-\zl_{\E}\op{id}_E)=\sum_{\E}Q_{\E}(\vec{X})(X_{\E}-\zl_{\E}\op{id}_E)$
is a polynomial in the $X_{\E}$. Up to factor $w,$ the {\small
LHS} reads
$$Q(\vec{X})(e)=p_1Q(\vec{X})i_1p_1(e)+p_2Q(\vec{X})i_1p_1(e)+p_1Q(\vec{X})i_2p_2(e)+p_2Q(\vec{X})i_2p_2(e),$$
where the second term of the {\small RHS} vanishes, in view of the
stability of $E_1$, and where the last term coincides with
$p_2(e),$ in view of Remark ($\star$) and Equation
(\ref{LeftInv2}). Eventually, cocycle $C=e\otimes w$ is
cohomologous to cocycle $$C_1=C-{\cal K}\zk C=\lp
p_1(e)-p_1Q(\vec{X})i_1p_1(e)-p_1Q(\vec{X})i_2p_2(e)\rp\otimes
w\in E_1\otimes\w.\mbox{  \rule{1.5mm}{2.5mm}}$$


The preceding proposition allows in particular recovering the
following well-known result:

\begin{cor}\label{dim0} Consider $n$ commuting endomorphisms $\vec{X}\in\op{End}^{\times
n}_{\C}(E)$ of a finite-dimensional complex vector space $E$, and
a point $\vec{\zl}\in\C^n$. Set
$\op{ker}\vec{X}_{\zl}:=\cap_{\E=1}^n\op{ker}(X_{\E}-\zl_{\E}\op{id})$.
If $\op{dim}(\op{ker}\vec{X}_{\zl})=0$, the Koszul cohomology
$KH^*(\vec{X}_{\zl},E)$ is trivial, and vice versa.
\end{cor}

{\it Proof}. It suffices to note that, due to Proposition
\ref{CharFinDimJointSpec}, the dimensional assumption means that
$\vec{\zl}\notin\zs(\vec{X})$, and to apply the preceding
proposition with $E_1=0$. Conversely, if there exists
$x\in\op{ker}\vec{X}_{\zl}\setminus\{0\}$, then ${\cal
K}_{\vec{X}_{\zl}}x=\sum_{\E=1}^n
(X_{\E}-\zl_{\E}\op{id})(x)\;\zh_{\E}=0$, so
that $x$ is a nonbounding $0$-cocycle. \rule{1.5mm}{2.5mm}\\

The next consequence of Proposition \ref{KoszCohoBasic} shows that
the Koszul cohomology $KH^*(\vec{X}_{\zl},E)$ is---roughly
spoken---made up by joint eigenvectors with eigenvalues
$\zl_{\E}$.\\

Consider $n$ commuting endomorphisms $\vec X=:\vec
X^{(1)}\in\op{End}^{\times n}_{\C}(E)$ of a finite-dimensional
complex vector space $E=:E^{(1)}=:F^{(1)}$, and a point
$\vec{\zl}\in\C^n$. For any $a\in\{2,3,\ldots\}$, if
$\op{ker}^{(a-1)}:=\op{ker}\vec X_{\zl}^{(a-1)}$ and
$E^{(a)}:=E^{(a-1)}/\op{ker}^{(a-1)}$, the
$$X_{\E}^{(a)}:=\left(
X_{\E}^{(a-1)}\right)^{\sharp},$$ $\E\in\{1,\dots,n\}$, defined
recursively by $X_{\E}^{(a)}=\left(
X_{\E}^{(a-1)}\right)^{\sharp}:E^{(a)}\ni [e^{(a-1)}]\raa
[X^{(a-1)}_{\E}e^{(a-1)}]\in E^{(a)}$, are again $n$ commuting
(well-defined) operators on a finite-dimensional complex vector
space. We iterate this procedure finitely many times, thus
obtaining operators $X_{\E}^{(a)}$, $a\in\{1,\ldots,s+1\}$, until
$\op{ker}^{(s+1)}=\op{ker}\vec X_{\zl}^{(s+1)}= 0$, or,
equivalently, $$\vec\zl\notin\zs(\vec X^{(s+1)}).$$ In the
following, we identify the operators $X_{\E}^{(a)}$ with their
models that arise from the choices of supplementary subspaces
$F^{(a)}$ of $\op{ker}^{(a-1)}$ in $E^{(a-1)}\simeq F^{(a-1)}$,
$a\in\{2,\ldots, s+1\}$, so that $E^{(a)}\simeq F^{(a)}\subset
E^{(a-1)}\simeq F^{(a-1)}$. If we denote by $i_a:F^{(a)}\raa
F^{(a-1)}$ the inclusion and by $p_a:F^{(a-1)}\raa F^{(a)}$ the
canonical projection, the isomorphism $E^{(a)}\simeq F^{(a)}$ is
$E^{(a)}\ni [f^{(a-1)}]\leftrightarrow p_af^{(a-1)}\in F^{(a)}$,
and operator $X_{\E}^{(a)}$, viewed as endomorphism of $F^{(a)}$,
reads \be X_{\E}^{(a)}=p_aX_{\E}^{(a-1)}i_a,\label{LowRightOp}\ee
since for any $f^{(a)}\in F^{(a)}$, we have
$X_{\E}^{(a)}f^{(a)}=X_{\E}^{(a)}[f^{(a)}]=[X_{\E}^{(a-1)}f^{(a)}]=p_aX_{\E}^{(a-1)}i_af^{(a)}$.

\begin{cor}\label{KoszCohoFin}
Let $\vec\zl\in\C^n$ be a point in $\C^n$, and let $\vec X=\vec
X^{(1)}\in\op{End}^{\times n}_{\C}(E)$ be $n$ commuting
endomorphisms of a finite-dimensional complex vector space
$E=F^{(1)}$. Denote by $\vec X^{(a)}\in\op{End}^{\times
n}_{\C}(F^{(a)})$, $a\in\{2,\ldots,s\}$, the above-depicted
``reduced'' operators on supplementary spaces $F^{(a)}$, and
denote by $\w=\w_n\la\vec{\zh}\ra$ the Grassmann algebra with $n$
generators $\vec{\zh}$.

Any cocycle
$$C\in E\otimes \w$$ of the Koszul complex $K^*(\vec{X}_{\zl},E)$
is cohomologous to a cocycle
$$C_1\in \lp
\op{ker}\vec{X}_{\zl}^{(1)}\oplus\op{ker}\vec{X}_{\zl}^{(2)}\oplus\ldots\oplus\op{ker}\vec{X}_{\zl}^{(s)}\rp\otimes\w.$$\end{cor}

{\it Proof}. It suffices to apply Proposition \ref{KoszCohoBasic}
to the obvious splitting
$$E=E_1\oplus E_2:=\left(\bigoplus_{a=1}^s\op{ker}^{(a)}\right)\oplus
F^{(s+1)}.$$ Indeed, the operators $\vec X_2$ considered in
Proposition \ref{KoszCohoBasic} read $X_{\E 2}=p_{s+1}\ldots
p_2X_{\E}i_2\ldots i_{s+1}=X_{\E}^{(s+1)}$, where we used the
afore-introduced notations $i_a$ and $p_a$. Hence, the spectral
condition $\vec\zl\notin\zs(\vec X_2)$ is satisfied by definition
of $s$, see above. Moreover, if $k^{(a)}\in \op{ker}^{(a)}\subset
F^{(a)}$, $a\in\{1,\ldots,s\}$, we have \be
X_{\E}\,k^{(a)}=X_{\E}\,i_2\ldots i_a\,k^{(a)}=p_{a}\ldots
p_2\,X_{\E}\,i_2\ldots
i_a\,k^{(a)}+\sum_{b=2}^{a}\zp_{b}\,p_{b-1}\ldots
p_2\,X_{\E}\,i_2\ldots i_a\,k^{(a)}.\label{Decomp}\ee Mapping
$\zp_b:F^{(b-1)}\raa \op{ker}^{(b-1)}$ is the second projection
associated with the decomposition $F^{(b-1)}=F^{(b)}\oplus
\op{ker}^{(b-1)}$, so that $\op{id}_{F^{(b-1)}}=p_b+\zp_b$. In
order to derive Equation (\ref{Decomp}), we utilized this upshot
for $b\in\{2,\ldots,a\}$. The first term of the {\small RHS} of
Equation (\ref{Decomp}) is
$X_{\E}^{(a)}k^{(a)}=\zl_{\E}\,k^{(a)}\in\op{ker}^{(a)}$, and the
terms characterized by index $b$ are elements of the spaces
$\op{ker}^{(b-1)}$.  Hence, space
$E_1=\oplus_{a=1}^s\op{ker}^{(a)}$ is stable under the action of
the $X_{\E}$ and Proposition \ref{KoszCohoBasic} is applicable.
\rule{1.5mm}{2.5mm}

\begin{cor} On the conditions of Corollary \ref{KoszCohoFin}, if for any $\E\in\{1,\ldots,n\}$, the kernel and the
image of the transformation $X_{\E}-\zl_{\E}\op{id}$ are
supplementary in $E$, then any cocycle $C\in E\otimes \w$ of the
Koszul complex $K^*(\vec{X}_{\zl},E)$ is cohomologous to a cocycle
$C_1\in \op{ker}\vec{X}_{\zl}\otimes\w.$\end{cor}

{\it Proof}. It suffices to prove that $s=1$. If $s\neq 1$, there
is a nonzero vector $x\in\op{ker}\vec X_{\zl}^{(2)}\subset
F^{(2)}$. Then, for any $k,\E\in\{1,\ldots,n\}$,
$(X_k-\zl_k\op{id})(X_{\E}-\zl_{\E}\op{id})x=(X_k-\zl_k\op{id})(p_2X_{\E}i_2x+\zp_2X_{\E}i_2x-\zl_{\E}x)
=(X_k-\zl_k\op{id})(\zp_2X_{\E}i_2x)=0,$ as
$\zp_2X_{\E}i_2x\in\op{ker}\vec X_{\zl}.$ Hence, for every $\E$,
we have $(X_{\E}-\zl_{\E}\op{id})x\in \op{ker}\vec
X_{\zl}\cap\op{im}(X_{\E}-\zl_{\E}\op{id})=0.$ Eventually,
$x\in(\op{ker}\vec X_{\zl})\cap F^{(2)}=0$, a contradiction.

\section{Koszul cohomology associated with Poisson cohomology}

We now come back to the Koszul cohomology implemented by a {\small
SRMI} tensor of $\R^n$. Let us recall that we deal with a {\small
SRMI} tensor $$\zL=\sum_{j<k}\alpha^{jk}Y_{jk}\quad
(\alpha^{jk}\in\R),$$ where the $Y_j$ are $n$ commuting linear
vector fields that verify $Y_{1\ldots n}\neq0$. The main building
block of the LP-cohomology of such a tensor has been identified as
the Koszul cohomology space $KH^*(\vec{X}_{\zd},E_r)$ associated
to the operators
$\vec{X}_{\zd}=(X_1-\zd_1\op{id},\ldots,X_n-\zd_n\op{id})$,
$X_j=\sum_{k}\alpha^{jk}Y_k$, $\za^{kj}=-\za^{jk}$,
$\zd_j=\op{div}X_j$ on the spaces $E_r={\cal S}^r\R^{n*}$,
$r\in\N$. We already pointed out that this cohomology can be
deduced from its complex counterpart
$KH^*(\vec{X}_{\zd}^{\C},E_r^{\C})$ (see Proposition
\ref{ComplexCohoKosz}), which is tightly related with joint
eigenvectors and the joint spectrum of $\vec{X}^{\C}$ or
$\vec{X}_{\zd}^{\C}$ (see Corollaries \ref{dim0} and
\ref{KoszCohoFin}). In this section, we further investigate the
Koszul cohomology space $KH^*(\vec{X}_{\zd}^{\C},E_r^{\C})$. In
particular, we reduce the computation of this central part
of the LP-cohomology space $LH^{*r}({\cal R},\zL)$ to essentially a problem of linear algebra,
and give a description of the spectrum of the transformations $\vec X_{\zd}^{\C}$.\\

When dealing with commuting operators on a finite-dimensional
complex vector space, it is natural to use an upper-triangular
representation of these transformations. The following theorem
shows that, for our endomorphisms $\vec X^{\C}_{\zd}$ of the space
$E_r^{\C}={\cal S}^r\C^{n*}$ (see below), which has the possibly
high (complex) dimension $N=(r+n-1)!/[r!(n-1)!]$ (if e.g. $r=10$
and $n=3$, this dimension equals $N=66$), the problem of finding
such a representation $\vec X^{\C}_{\zd}\in\op{gl}(N,\C)^{\times
n}$ (we denote the operators and their representation by the same
symbol) reduces to the quest for an upper-triangular
representation $\vec a\in\op{gl}(n,\C)^{\times n}$ of some
commuting transformations $\vec a$ of $\C^n$. More precisely, the
$a_j$, $j\in\{1,\ldots,n\}$, are the commuting matrices
$a_j=(J^1)^{-1}Y_j\in\op{gl}(n,\R)$ that correspond to the
commuting linear vector fields $Y_j.$

\begin{prop}\label{IndBasis} Any basis of $\C^n$, in which the commuting operators $\vec
a$ have an upper-triangular representation, naturally induces a
basis of $E_r^{\C}={\cal S}^r\C^{n*}$, in which all the
transformations $\vec X_{\zd}^{\C}$ are
upper-triangular.\end{prop}

Let us first mention that in the sequel the use of super- and
subscripts is dictated by esthetic criteria and not at all by
contra- or covariance.\\

{\it Proof}. In the following, we denote by $x=(x_1,\ldots,x_n)$
(resp. $z=(z_1,\ldots,z_n)$) the points of $\R^n$ (resp. $\C^n$)
as well as their coordinates in the canonical basis
$(e_1,\ldots,e_n)$. As usual, we set
$Y_k=\sum_m\E_{km}\p_{x_m}=\sum_{mp}a_{k}^{mp}\,x_p\p_{x_m}$ and
use notations as $x^{\zb}=x_1^{\zb_1}\ldots x_n^{\zb_n}$,
$\zb\in\N^n$.

The complexification $E_r^{\C}$ of
$$E_r={\cal S}^r\R^{n*}=\{P\in\Ci(\R^n):P(x)=\sum_{\vert\zb\vert=r}r_{\zb}x^{\zb}\quad (x\in\R^n,
r_{\zb}\in\R)\}$$ is
$$E_r\oplus i E_r\simeq E_r^{\C}\simeq {\cal S}^r\C^{n*}=\{P\in\Ci(\C^n):P(z)=\sum_{\vert\zb\vert=r}c_{\zb}z^{\zb}\quad
(z\in\C^n, c_{\zb}\in\C)\}.$$ It is also easily seen that the
complexification $Y_k^{\C}\in\op{End}_{\C}(E_r^{\C})$ of
$Y_k\in\op{End}_{\R}(E_r)$ is the holomorphic vector field
$$Y_k^{\C}=\sum_{mp}a_{k}^{mp}\,z_p\p_{z_m}\in\op{Vect}^{10}(\C^n)$$ of $\C^n$.

It is well-known that the $n$ commuting matrices
$a_j=(J^1)^{-1}Y_j\in\op{gl}(n,\R)$ can be reduced simultaneously
to upper-triangular matrices by a unitary matrix
$U\in\op{U}(n,\C).$ Consider any matrix $U\in\op{GL}(n,\C)$ (resp.
any basis $(e'_1,\ldots,e'_n)$ of $\C^n$), such that the
$b_j=U^{-1}a_jU\in\op{gl}(n,\C)$ are upper-triangular (resp. in
which the transformations $\vec a$ are all upper-triangular).
Denote by ${\frak z}=({\frak z}_1,\ldots,{\frak z}_n)$ the
components of the vectors $z=\sum_j{\frak z}_je'_j\in\C^n$ in the
basis $(e_1',\ldots,e_n')$, and let $(\ze_1',\dots,\ze'_n)$ be the
dual basis of this new basis. If viewed as a basis of the space
$E_r^{\C}$ of degree $r$ homogeneous polynomials of $\C^n$, the
induced basis $\ze'_{j_1}\vee\ldots\vee\ze'_{j_r}$,
$j_1\le\ldots\le j_r,$ of the space ${\cal S}^r\C^{n*}$ of
symmetric covariant $r$-tensors of $\C^n$ reads ${\frak z}^{\zb},$
$\zb\in\N^n,$ $\vert\zb\vert=r$.

In order to find the matrices of the operators $\vec X_{\zd}^{\C}$
in this ``natural'' basis ${\frak z}^{\zb},$ $\zb\in\N^n,$
$\vert\zb\vert=r$ of $E^{\C}_r$, we range the vectors ${\frak
z}^{\zb}$ according to the lexicographic order $\prec$ and perform
the coordinate change $z=U{\frak z}, \p_{z}=\widetilde{\p_{\frak
z}z}^{-1}\p_{\frak z}$ in the first order linear differential
operators $(X_j-\zd_j\op{id})^{\C}$. We get
\begin{eqnarray*}(X_j-\zd_j\op{id})^{\C}&=&\sum_k\za^{jk}\sum_{m\le
p}b_{k}^{mp}\,{\frak z}_p\p_{{\frak z}_m}-\zd_j\op{id}^{\C}\\
&=&\sum_{km}\za^{jk}b_{k}^{mm}\,\left({\frak z}_m\p_{{\frak
z}_m}-\op{id}^{\C}\right)+\sum_k\sum_{m<p}\za^{jk}b_{k}^{mp}\,{\frak
z}_p\p_{{\frak z}_m},
\end{eqnarray*} since
$\zd_j=\op{div}X_j=\sum_{km}\za^{jk}a_k^{mm}=\sum_{km}\za^{jk}b_k^{mm}.$
As the image of vector ${\frak z}^{\zb}$ by operator
$(X_j-\zd_j\op{id})^{\C}$ is \be(X_j-\zd_j\op{id})^{\C}{\frak
z}^{\zb}=\sum_{km}\za^{jk}b_{k}^{mm}\,\left(\zb_m-1\right){\frak
z}^{\zb}+\sum_k\sum_{m<p}\za^{jk}b_{k}^{mp}\zb_m\,{\frak
z}^{\zb-e_m+e_p},\label{Spect}\ee where ${\frak
z}^{\zb-e_m+e_p}\prec{\frak z}^{\zb}$, the matrices of the
commuting operators $(X_j-\zd_j\op{id})^{\C}$,
$j\in\{1,\ldots,n\}$, in the basis ${\frak z}^{\zb}$,
$\zb\in\N^n,$ $\vert\zb\vert=r$, of space
$E_r^{\C}$, are all upper-triangular. \rule{1.5mm}{2.5mm}\\

The next theorem provides a description of the joint spectrum
$\zs_r(\vec{X}_{\zd}^{\C})$ of the operators $\vec X_{\zd}^{\C}\in\op{End}_{\C}^{\times n}(E_r^{\C})$.\\

Let $B\in\op{gl}(n,\C)$ be the matrix $B_{jk}=b_{j}^{kk}$ made up
by the diagonals of the matrices $b_j,$ see above.

\begin{theo}\label{JointSpecPoiss} The joint spectrum $\zs_r(\vec{X}_{\zd}^{\C})$ of the commuting operators
$\vec{X}_{\zd}^{\C}\in\op{End}^{\times n}_{\C}(E_r^{\C})$ on the
finite-dimensional complex vector space $E_r^{\C}$, is given by
$$\zs_r(\vec{X}_{\zd}^{\C})=\{\za BI:I\in(\N\cup\{-1\})^n,\vert I\vert=r-n\}\subset\C^n,$$
where $\vert I\vert=\sum_j I_j$ denotes the length of $I$.
\end{theo}

{\it Proof}. Direct consequence of Proposition
\ref{CharFinDimJointSpec} and Equation (\ref{Spect}).
\rule{1.5mm}{2.5mm}\\

\noindent {\bf Remark}. In Proposition \ref{JointEigen}, we showed
that for all $k$, $Y_kD=(\op{div}Y_k)D$, where $D=\op{det}\E\in
E_n\subset E_n^{\C}$. It of course follows that for all $j$,
$X_j^{\C}D=X_jD=(\op{div}X_j)D=\zd_j\op{id}^{\C}D,$ so that
$\vec{0}=(0,\ldots,0)\in\zs_n(\vec{X}_{\zd}^{\C})$. This last
upshot is immediately recovered from Theorem
\ref{JointSpecPoiss}.\\

Set $K_r(\vec{X}_{\zd}^{\C})=\{I\in\op{ker}(\alpha
B):I\in\left(\N\cup\{-1\}\right)^n,\vert I\vert=r-n\}$. Corollary
\ref{dim0} can then be reformulated as follows.

\begin{cor}\label{dim0'}
The Koszul cohomology $KH^*(\vec{X}_{\zd}^{\C},E_r^{\C})$ is
acyclic if and only if $K_r(\vec{X}_{\zd}^{\C})=\varnothing$.
\end{cor}

{\it Proof}. Indeed, $KH^*(\vec{X}_{\zd}^{\C},E_r^{\C})$ is
trivial if and only if $\op{dim}(\op{ker}\vec{X}_{\zd}^{\C})=0$,
if and only if $\vec{0}\notin\zs_r(\vec{X}_{\zd}^{\C})$, i.e. if
and only if $K_r(\vec{X}_{\zd}^{\C})=\varnothing$.
\rule{1.5mm}{2.5mm}\\

We now depict a convenient method that allows finding a basis of
the space
$$\op{ker}\vec{X}_{\zd}^{\C\,(1)}\oplus\op{ker}\vec{X}_{\zd}^{\C\,(2)}
\oplus\ldots\oplus\op{ker}\vec{X}_{\zd}^{\C\,(s)},$$ which houses
the Koszul cohomology $KH^*(\vec{X}_{\zd}^{\C},E_r^{\C})$, see
Corollary \ref{KoszCohoFin}.\\

In order to simplify notations, we systematically omit in the
following description superscript $\C$. We write e.g. $\vec
X_{\zd},E_r,\ldots$ instead of $\vec
X_{\zd}^{\C},E_r^{\C},\ldots$\\

Consider any basis $(e_1,\ldots,e_N)$ of $E_r$ that generates an
upper-triangular representation $T_1,\ldots,T_n$ of the operators
$\vec X_{\zd}$. The kernel $\op{ker}\vec X_{\zd}$ is then
described by the $n$ triangular systems \be
T_1Z=0,\ldots,T_nZ=0,\label{BasisKer}\ee (each one) of $N$
homogeneous linear equations in the $N$ complex unknowns
$Z=(Z^1,\ldots,Z^N)$.

As understood before, $\vec 0\in\zs_r(\vec X_{\zd})$ if and only
if at least one of the lines
$\vec{T}\,^q=(T_1^{qq},\ldots,T_n^{qq})$, $q\in\{1,\ldots,N\}$, is
$\vec 0=(0,\ldots,0)$. We refer to the number $\zm$ of such $\vec
0$--lines $\vec T\,^{q_1},\ldots,\vec T\,^{q_{\zm}}$, $q_1<\ldots<
q_{\zm}$, as the multiplicity of $\vec 0$ in the spectrum
$\zs_r(\vec X_{\zd})$ (in the considered basis
$(e_1,\ldots,e_N)$). Of course, the general solution of System
(\ref{BasisKer}) is a linear combination $Z=\sum_jc_jK_j$,
$c_j\in\C$, of $d=\op{dim}\op{ker}\vec X_{\zd}$ independent
vectors $K_j\in\C^N$. Let \be k_j=K_j^1e_1+\ldots
+K_j^{q_{\zn_j}}e_{q_{\zn_j}},\;j\in\{1,\ldots,
d\},\label{BasisKerVect}\ee be the corresponding basis of
$\op{ker}\vec X_{\zd}$. It can quite easily be seen---just
``solve'' System (\ref{BasisKer}) and start imagining a
configuration that leads to the maximal dimension of the space of
solutions---that $d\le \zm$ and that the components
$K_j^{q_{\nu_j}}\neq 0$ of the vectors $k_j$ with highest
superscript correspond to $\vec 0$--lines
${q_{\nu_1}}<\ldots<{q_{\nu_d}}$.

The $N$-tuple
$(k_1,\ldots,k_d,e_1,\ldots,\widehat{e_{q_{\nu_1}}},\ldots,\widehat{e_{q_{\nu_d}}},\ldots,e_N)$
is a basis of $E_r$, since the determinant in the basis
$(e_1,\ldots,e_N)$ of the permuted $N$-tuple
$(e_1,\ldots,k_1,\ldots,k_d,\ldots,e_N)$ equals
$K_1^{q_{\nu_1}}\ldots K_d^{q_{\nu_d}}\neq 0$. Observe that the
$k_j$ are joint eigenvectors of the $\vec X_{\zd}$ associated with
eigenvalue $0$. Moreover, in view of Equation
(\ref{BasisKerVect}), every vector $e_{q_{\nu_j}}$ can be written
in terms of ``lower'' vectors of the new basis. Hence, the first
$d$ columns of the representative matrices $T'_1,\ldots,T'_n$ of
the operators $\vec X_{\zd}$ in the new basis vanish, these
matrices are again upper-triangular, and the lines $\vec T^{q}$,
$q\in\{1,\ldots,N\}$, are unchanged up to permutation. The
matrices $T'_{\E}+\zd_{\E}\op{id}\in\op{gl}(N,\C)$ correspond to
the operators $X_{\E}$, $\E\in\{1,\ldots,n\}$, and their lower
right submatrices
$(T_{\E}^{'}+\zd_{\E}\op{id})^{(2)}\in\op{gl}(N-d,\C)$ (resp.
$T_{\E}'^{(2)}$) correspond to the operators $X_{\E}^{(2)}$ (resp.
$X_{\E}^{(2)}-\zd_{\E}\op{id}^{(2)}$), see Equation
(\ref{LowRightOp}) and Corollary \ref{KoszCohoFin}.

In other words, in the basis
$(e_1,\ldots,\widehat{e_{q_{\nu_1}}},\ldots,\widehat{e_{q_{\nu_d}}},\ldots,e_N)$
of a space $F^{(2)}_r$, see Corollary \ref{KoszCohoFin}, which is
supplementary to $\op{ker}\vec X_{\zd}$ in $E_r$, the operators
$\vec X_{\zd}^{(2)}$ are represented by upper-triangular matrices
$T_1'^{(2)},\ldots, T_n'^{(2)}$. Thus, the above-detailed
procedure can be iterated and the general solution of another
packet of $n$ (smaller) triangular systems of linear equations \be
T_1'^{(2)}Z=0,\ldots, T_n'^{(2)}Z=0,\label{BasisKer2}\ee provides
a basis $k_1^{(2)},\ldots,k_{d_2}^{(2)}$ of $\op{ker}\vec
X_{\zd}^{(2)}$, et cetera.\\

\noindent{\bf Remarks.} \begin{itemize}\item {\it The solutions of
the triangular systems of homogeneous linear equations
(\ref{BasisKer}), (\ref{BasisKer2}), ... generate a basis of the
locus
$$\op{ker}\vec{X}_{\zd}^{\C\,(1)}\oplus\op{ker}\vec{X}_{\zd}^{\C\,(2)}
\oplus\ldots\oplus\op{ker}\vec{X}_{\zd}^{\C\,(s)}$$ of the Koszul
cohomology space $KH^*(\vec X_{\zd}^{\C},E_r^{\C})$}. \item
Observe that if the $b_{\E}=U^{-1}a_{\E}U$ have been computed, the
upper-triangular matrix representations $T_1,\ldots,T_n$ of the
transformations $\vec X_{\zd}^{\C}$ in the corresponding basis
${\frak z}^{\zb}$, $\zb\in\N^n,$ $\vert\zb\vert=r$, of $E_r^{\C}$
are known, see Equation (\ref{Spect}), and explicit computations
can actually be performed.\item As the multiplicity of $\vec 0$ in
the spectrum of the endomorphisms $\vec X_{\zd}^{\C\,(2)}$ is
$\zm-d$, and as its multiplicity in the spectrum of the $\vec
X_{\zd}^{\C\,(s+1)}$ vanishes, by definition of $s$, we get
\be\zm=d+d_2+\ldots+d_s=\sum_{j=1}^{s}\op{dim}\op{ker}\vec
X_{\zd}^{\C\,(j)},\label{MultiplicDim}\ee with self-explaining
notations. As the {\small RHS} of this equation is independent of
the considered basis, the concept of multiplicity of a point
$\zl\in\C^n$ in the joint spectrum of commuting transformations of
a finite-dimensional vector space, makes sense. Although this
result might be well-known, we could not find it anywhere in
literature.
\end{itemize}

\noindent {\bf Example 1.} Consider structure $\zL_2$ of the
{\small DHC}, see Theorem \ref{ClassTheo}, and assume that $a\neq
0,b=0$. It is easily checked that the matrix $$ U=\left(
  \begin{array}{ccc}
    0 & \frac{i}{\sqrt{2}} & \frac{-i}{\sqrt{2}} \\
    0 & \frac{1}{\sqrt{2}} & \frac{1}{\sqrt{2}} \\
    1 & 0 & 0 \\
  \end{array}
\right)$$ transforms the above-mentioned matrices $a_{\E}$
simultaneously into upper-triangular matrices $b_{\E}$. 
A short computation yields that the space
$K_{3t}(\vec{X}_{\zd}^{\C})$, $t\in\N$, contains the unique point
$I_t=(t-1,t-1,t-1)$, so that the multiplicity $\zm$ of $\vec 0$ in
the joint spectrum $\zs_{3t}(\vec X_{\zd}^{\C})$ equals $1$, see
proof of Theorem \ref{JointSpecPoiss}. It follows that the Koszul
cohomology spaces $KH^*(\vec{X}_{\delta}^{\C},E_{3t}^{\C})$ are
not trivial, see Corollary \ref{dim0'}. Furthermore, since the
matrices $b_{\E}$ are in fact diagonal in this example, Equation
(\ref{Spect}) entails that ${\frak z}_1^t{\frak z}_2^t{\frak
z}_3^t$ belongs to the kernel $\op{ker}_{3t}\vec{X}_{\zd}^{\C}$ of
operators $\vec X_{\zd}^{\C}$ in space $E_{3t}^{\C}$. If we take
into account Equation (\ref{MultiplicDim}), we see that
$\op{ker}_{3t}\vec{X}_{\delta}^{\C}=\C {\frak z}_1^t{\frak
z}_2^t{\frak z}_3^t$ and that the reduced operators $\vec
X_{\zd}^{\C\,(j)}$, $j\in\{2,\ldots,s\}$, do not exist, i.e. that
$s=1$. Hence, and since the change to canonical coordinates is
$z=U{\frak z}$, see proof of Proposition \ref{IndBasis}, the
cohomology space $KH^p(\vec{X}_{\delta}^{\C},E_{3t}^{\C})$,
$p\in\{0,1,2,3\}$, $t\in\N,$ is located inside
$${\frak z}_1^t{\frak z}_2^t{\frak z}_3^t\bigoplus_{j_1<\ldots<j_p}\C Y_{j_1\ldots
j_p}=(z_1^2+z_2^2)^tz_3^t\bigoplus_{j_1<\ldots<j_p}\C Y_{j_1\ldots
j_p}.$$ This rather easily obtained upshot is in accordance with
the results of \cite{MP} (modulo slight changes in definitions and notations [e.g. the roles of parameters $a$ and $b$ are exchanged]).\\

\noindent{\bf Example 2.} For structure $\zL_3$ of the {\small
DHC} and parameter value $a=0,$ depending on the value of $r$, the
multiplicity of $\vec 0$ in the spectrum $\zs_r(\vec
X_{\zd}^{\C})$ equals $0$ or $1$---and computations are similar to
those of the preceding example---, except in the case $r=3$, which
generates multiplicity $3$. Since for $\zL_3$ the matrices
$a_{\E}$ are lower-triangular, a coordinate change
$z\leftrightarrow {\frak z}$ is not necessary and it can easily be
seen that we have $s=3$ and $$\op{ker}_3\vec X_{\zd}^{\C}=\C
z_1^2z_3,\op{ker}_3\vec X_{\zd}^{\C\,(2)}=\C
z_1z_2z_3,\op{ker}_3\vec X_{\zd}^{\C\,(3)}=\C z_2^2z_3.$$ The
corresponding cohomological upshots are part of the computation of
the LP-cohomology of $\zL_3$ that we detail in the next section.\\

\noindent {\bf Remark.} Remember that the operators $X_i$ are
defined by $X_i=\sum_{j}\za^{ij}Y_j$, with $\za^{ji}=-\za^{ij}$.
Hence, matrix $\za\in\op{gl}(n,\R)$ is skew-symmetric, and
$\op{det}\za$ vanishes for odd $n$. Of course, the corresponding
non-trivial linear combination $\sum_ic_i\za^{i*}=0$, induces a
non-trivial combination $\sum_ic_iX_i=0$ of the $X_i$ (and the
$X_i-\zd_i\op{id}$), which is significant in computations. In the
even dimensional ($n=2m,m\in\{2,3,\ldots\}$) maximal rank
($\op{rk}\za=n$) case, the Koszul cohomology $KH^*(\vec
X_{\zd}^{\C},E_r^{\C})$ has the following simple description. If
$($in even dimension $n$$)$ $\op{det}\alpha\neq 0$, then
$$\bigoplus_{r\in\N}KH^{0}(\vec X_{\zd}^{\C},E_r^{\C})=\C\,{\cal D},$$
where ${\cal D}$ denotes the complex clone of $\op{det}\E$, and,
for any $r\neq n$ and any $p\in\{1,\ldots,n\}$, the cohomology
space $KH^{p}(\vec X_{\zd}^{\C},E_r^{\C})$ vanishes. We do not
detail the proof that is, roughly, along the lines of Proposition
\ref{JointEigen}. If $r=n$, the situation is more complicated and
new elements of
$\op{ker}\vec{X}_{\zd}^{\C\,(1)}\oplus\op{ker}\vec{X}_{\zd}^{\C\,(2)}
\oplus\ldots\oplus\op{ker}\vec{X}_{\zd}^{\C\,(s)}$ may enter the
play. \label{PoissKoszSect}

\section{Cohomology spaces of structures $\zL_3$ and $\zL_9$}

We already pointed out that the LP-cohomology (or ${\cal
R}$-cohomology) of {\small SRMI} tensors can be deduced from a
Koszul cohomology (${\cal P}$-cohomology) and a relative
cohomology (${\cal S}$-cohomology), see Theorem
\ref{PoissKoszRel}, Theorem \ref{PoissKosz}, and Proposition
\ref{SuppRelat}.

The involved Koszul cohomology has been studied in the last
section. We particularized our upshots by means of (pertinent)
examples, see Examples 1 and 2, Section \ref{PoissKoszSect}.

Within the cohomology computations of {\small SRMI} tensors of the
{\small DHC}, ${\cal S}$-cohomology has so far been determined
``by hand''. In the majority of cases, the LP-cohomology operator
respects, in addition to the degrees $p$ and $r$, a partial
polynomial degree $k$ (e.g. the coboundary operator associated
with $\zL_3$ respects the partial degree in $x=x_1,y=x_2$), so
that we can decompose space ${\cal S}^{pr}$ into smaller spaces
${\cal S}_{kr}^{p}$ (made up by the elements of ${\cal S}^{pr}$
that have partial degree $k$), see \cite{MP}. The cohomology
operator of structure $\zL_9$ however, does not respect any
additional degree. The ${\cal S}$-cohomology of $\zL_9$ is
therefore quite intricate.

Theorem \ref{PoissKoszRel} leads to the following cohomological
upshots for structures $\zL_3$ and $\zL_9$. No proofs will be
given (for a description of an application of the technique, see
\cite{MP}).

\begin{theo}\label{PoissCohoL3}
If $a\neq 0$, the cohomology spaces of structure $\zL_3$ are
$$LH^{0*}({\cal R},\zL_3)=\R,$$
$$LH^{1*}({\cal R},\zL_3)=\R Y_1+\R Y_2+\R Y_3,$$
$$LH^{2*}({\cal R},\zL_3)= \mathbb{R}Y_{23} \oplus \mathbb{R}Y_{31}\oplus \mathbb{R} (2yz\p_{31} +
y^2\p_{12}),$$
$$LH^{3*}({\cal R},\zL_3)=\mathbb{R} \, \p_{123} \oplus \mathbb{R} \, y^2z \, \p_{123},$$
where the $Y_i$ are those defined in Theorem \ref{ClassTheo}.
\end{theo}


\begin{theo}\label{PoissCohoL9}
If $a\neq0$, the cohomology spaces of structure $\zL_9$ are
$$LH^{0*}({\cal R},\zL_9)=\R,$$
$$LH^{1*}({\cal R},\zL_9)=\R Y_1+\R Y_2+\R Y_3,$$
$$LH^{3*}({\cal R},\zL_9)=\oplus_{r\in\N}\R z^r\p_{123},$$ and
$$LH^{2*}({\cal R},\zL_9)=\oplus_{r\in\N}H_r^2,$$ where
\begin{eqnarray*}
  &&H_0^2=\R \p_{23},\quad H_1^2=\R C_1^0,\quad H_3^2=\R C_1^2,  \\
  &&H_2^2=\R x^{2}\p_{23}+\R xz(\p_{23}-2^{-1}\p_{31})+\R(xz\p_{12}-z^2\p_{23})\\
  &&\quad\quad\quad+\R(yz\p_{12}+(-27a^2x^2-9axz+3ay^2-z^2)\p_{31}),\\
  &&H^2_{r+1}=\R C_1^r+\R C_2^r,\quad r\geq 3,
\end{eqnarray*}
with
\begin{eqnarray*}C_1^r&=&-a(xz^r+ry^2z^{r-1})\p_{12}
+(9a^2xy^r+a(3r-1)(r+1)^{-1}z^{r+1})\p_{23}\\
&&+ayz^r\p_{31}
\end{eqnarray*}
and
\begin{eqnarray*}C_2^r&=&(-a(r-2)y^4z^{r-3}+y^2z^{r-1})\p_{12}\\&&
+(9a^2xy^2z^{r-2}-9ar^{-1}xz^r
+3a(r-3)(r-1)^{-1}y^2z^{r-1}-3(r-1)r^{-1}(r+1)^{-1}z^{r+1})\p_{23}\\&&
+(6a(r-1)^{-1}xyz^{r-1}-ay^3z^{r-2}-r^{-1}yz^r)\p_{31},
\end{eqnarray*}
where the $Y_i$ are those defined in Theorem \ref{ClassTheo} $($
and where the terms that contain a power of $x,$ $y,$ or $z$ with
a negative exponent are ignored $)$.
\end{theo}

\section{Cohomological phenomena}

Let us outline the most important cohomological phenomena.\\

Consider a {\small SRMI} Poisson structure $\zL=\sum_{i<j}\za^{ij}Y_{ij}$.\\

It is easily checked that the curl vector field of $\zL$, see
Section \ref{admstr}, is given by $K(\zL)=\sum_i\zd_iY_i$,
$\zd_i=\op{div}X_i$, $X_i=\sum_j\za^{ij}Y_j$. Consequently,
K-exactness is (in $\R^n$, $n\ge 3$) equivalent with
divergence-freeness. Note now that the $0$--cohomology space
$LH^{0*}({\cal R},\zL)$ of $\zL$, or space $\op{Cas}(\zL)$ of
Casimirs of $\zL$, coincides with the kernel $\op{ker}\vec X$, see
Equation \ref{SRMICobFin}. Hence, in view of Proposition
\ref{JointEigen}, for a K-exact tensor, $D^p=(\op{det}\E)^p$ is a
joint eigenvector of the $X_i$ with eigenvalues $p\,\zd_i=0$, i.e.
$D^p\in\op{ker}\vec X$. It follows that, for K-exact {\small SRMI}
Poisson tensors,
$$\oplus_{p\in\N}\R D^p\subset LH^{0*}({\cal
R},\zL)=\op{Cas}(\zL).$$

As for the $1$--cohomology space $LH^{1*}({\cal R},\zL)$, let us
first remark that the stabilizer ${\frak g}_{\zL}$, viewed as a
Lie subalgebra of linear vector fields ${\cal X}^1_0(\R^n)$, is
made up by $1$-cocycles (by definition) that do not bound (degree
argument), i.e.
$${\frak g}_{\zL}\subset LH^{1*}({\cal R},\zL).$$

Moreover, as LP-cohomology is an associative graded commutative
algebra, the classes of the cocycles in
$$\op{Cas}(\zL)\otimes\w^p\frak{g}_{\zL},$$
$0\leq p\leq n$, are ``preferential'' LP-cohomology classes. Such
classes massively appear in the LP-cohomology of {\small SRMI}
tensors of the {\small DHC}, see \cite{MP}, and of twisted {\small
SRMI} tensors, see \cite{AMPN}.\\

However, two other types of classes systematically appear in
LP-cohomology.
\begin{enumerate}\item The classes of type I originate from ${\cal P}$-cohomology.
In fact, roughly spoken, the locus
$\op{ker}\vec{X}_{\zd}^{\C\,(1)}\oplus\op{ker}\vec{X}_{\zd}^{\C\,(2)}
\oplus\ldots\oplus\op{ker}\vec{X}_{\zd}^{\C\,(s)}$ of the Koszul
cohomology associated with the considered Poisson cohomology
generates in some cases nonbounding cocycles in ${\cal
R}$-cohomology. For instance, for structure $\zL_7$, the rational
functions $D'^{\frac{\gamma}{2}}z^{-1}$, $D'=x^2+y^2$, $\gamma\in
2\N^*$, induce the classes $D'^{\frac{\gamma}{2}}z^{-1}Y_3$,
$Y_3=z\p_3$, in space $LH^{1*}({\cal R},\zL_7)$. \item The classes
of type II are due to ${\cal S}$-cohomology. Indeed, let $\frak s$
be a cochain in space ${\cal S}$, which is supplementary to ${\cal
R}$ in ${\cal P}$. It happens that $\p_{\zL}\frak s\in{\cal R}$.
Then, $\p_{\zL}{\frak s}$---a coboundary of a cochain from the
outside of ${\cal R}$---is typically a nonbounding cocycle in
${\cal R}$.\end{enumerate}

We refer to these two types of cohomology classes as ``singular
classes'', since some of their coefficients are polynomials on the
singular locus of the considered Poisson tensor.\\

Let us finally briefly comment on the impact of LP- and
K-exactness on the structure of LP-cohomology. If tensor $\zL$, or
part of this tensor, is LP-exact, see Section \ref{admstr}, some
elements of space $\w^2\frak{g}_{\zL}$ may be bounding cocycles.
For instance, part $Y_{12}$ of structure $\zL_3$ of the {\small
DHC} is LP-exact and disappears in the second cohomology space,
see Theorem \ref{PoissCohoL3}. Hence, LP-exactness impoverishes
LP-cohomology. In view of the above remark on Casimir functions
and the observations made in earlier works, we know that
K-exactness significantly enriches the cohomology. Therefore,
richness of LP-cohomology depends in some sense on the distance of
the Poisson tensor to LP- and K-exactness.

\end{document}